\documentclass[12pt]{amsart}
\usepackage{amsmath,amsfonts,amscd,amssymb,latexsym,amsthm}
\usepackage{amsmath,amsfonts,amscd,amssymb,latexsym}
\title[Free Partially Commutative
Groups]{Divisibility Theory
and Complexity of Algorithms for Free Partially Commutative
Groups}
\author[E. S. Esyp]{Evgenii S. Esyp}

\address{E. S. Esyp, I. V. Kazachkov,
V. N. Remeslennikov: Omsk Branch of Mathematical Institute SB RAS,
13 Pevtsova Street,  Omsk 644099,  Russia}
\author[I. V. Kazachkov]{Ilya V. Kazachkov}
\thanks{The second author was partially supported by the
London Mathematical Society and by EPSRC grant GR/S71200.}

\author[V. N. Remeslennikov]{Vladimir N. Remeslennikov}\thanks{The third author was
supported by EPSRC grants GR/R29451 and GR/S71200.}

\newtheorem{lem}{Lemma}[section]
\newtheorem{thm}[lem]{Theorem}

\newtheorem{prop}[lem]{Proposition}
\newtheorem*{propn}{Proposition}

\theoremstyle{definition}
\newtheorem{alg}[lem]{Algorithm}

\newtheorem{defn}[lem]{Definition}
\newtheorem*{defnn}{Definition}
\newtheorem{expl}[lem]{Example}
\newtheorem{rem}[lem]{Remark}

\newcommand{\1}{\ensuremath{ {}^\circ}}
\newcommand{\inv}{\ensuremath{{}^{-1}}}
\newcommand{\Gplus}{\ensuremath{{\mathbb{G}}^+}}

\newcommand{\bi}{\begin{itemize}}
\newcommand{\ei}{\end{itemize}}
\newcommand{\bq}{\begin{quote}}
\newcommand{\eq}{\end{quote}}
\newcommand{\beq}{\begin{equation}}
\newcommand{\eeq}{\end{equation}}
\newcommand{\bea}{\begin{eqnarray*}}
\newcommand{\eea}{\end{eqnarray*}}

\newcommand{\ad}{\mathop{{\rm ad}}}
\newcommand{\az}{\mathop{{\ensuremath{\alpha}}}}
\newcommand{\cred}{\mathop{{\rm cr}}}
\newcommand{\CR}{\mathop{{\rm CR}}}
\newcommand{\gd}{\mathop{{\rm gd}}}
\newcommand{\lm}{\mathop{{\rm lm}}}
\newcommand{\md}{\mathop{{\rm gd}}}
\newcommand{\n}{\mathop{{\rm n}}}
\newcommand{\rel}{\mathop{{\rm rel}}}

\newcommand{\rn}{\mathop{{\rm rn}}}
\newcommand{\sign}{\mathop{{\rm sign}}}

\newcommand{\cC}{\ensuremath{\mathcal{C}}}

\newcommand{\pP}{\ensuremath{\mathcal{P}}}

\newcommand{\BA}{\ensuremath{\mathbb{A}}}
\newcommand{\GG}{\ensuremath{\mathbb{G}}}
\newcommand{\MM}{\ensuremath{\mathbb{M}}}

\newcommand{\ZZ}{\ensuremath{\mathbb{Z}}}

\begin{document}
\begin{abstract}
The original version of the paper was published in Contemporary
Mathematics 378 ``Groups, Languages, Algorithms''; 2005, pp.
319-348. This is a modified version with Appendix that holds a
corrected formulation of Proposition 4.1.
\end{abstract}

\maketitle

\tableofcontents

\section{Introduction}

Free partially commutative groups  (or \emph{partially commutative
groups}, for brevity) have many remarkable properties.

\begin{itemize}
    \item They arise naturally in many branches of mathematics and computer
    science. This led to a variety of names under which they are
    known: \emph{semifree groups} \cite{Baudisch2,Baudisch}, \emph{graph groups} \cite{Dr1,surf,Serv,VanWyk},
    \emph{right-angled Artin groups} \cite{BeBr,CW}, \emph{locally free groups}
    \cite{DeN,Vershik}.
    \item  Normal forms of elements in these groups are  extremely
    convenient from computational point
of view.
    \item  All major algorithmic problems for partially commutative groups are
solvable.
    \item The groups have very rich subgroup structure.
    In particular, the fundamental group of almost every
surface is a subgroup of a suitable partially commutative group
(see \cite{CW}).
    \item Free partially commutative  groups have interesting statistical
properties, and many of their statistical parameters can be
estimated.
    \item  Right-angled   groups are intrinsically connected to other
types of  Artin groups \cite{collins,pride}.
\end{itemize}

Not surprisingly, the rich pickings of the theory of partially
commutative groups attract researchers from various areas of
mathematics, and the group are being investigated from a
refreshing variety of points of view.

Without attempting to give a full survey of the existing theory,
we mention that the structural theory of partially commutative
groups was developed in \cite{DK,Serv} and structure of subgroups
studied in \cite{CW,surf}. Paper \cite{BeBr} discusses
applications to geometry and topology, while major statistical
characteristics for these groups are estimated in \cite{DeN,
Vershik, Vershik1} and the Poisson-Furstenberg boundary described
in \cite{Malyutin}. Some algorithmic problems for right angled
groups were solved in \cite{DTh,Wrath}.

We can summarize the main results of this paper sa follows. Let
$G$ be an arbitrary partially commutative group.
\begin{itemize}

\item We develop divisibility theory for $G$
(and not just for semigroups of positive elements, as in the
classical paper on Artin groups \cite{BrSa}). It closely resembles
the classical divisibility theory of integers. In particular, we
have reasonable concepts of the least common multiple and greatest
common divisor.

    \item We prove the existence of polynomial algorithms for solving
    main algorithmic problems and estimate the boundaries for their
    complexity.

    \item With the help of the algorithms of divisibility theory
    in each class $w^{G}$ of all elements conjugate to
    the element $w$ we compute the canonical representative.
    Its normal formal is called in the paper the \emph{minimal exhausted form} of $w^G$.
    This form is easy from the
    computational point of view (see Outline of this algorithm in Section
    \ref{ss:conrelsub}).

    \item These results allow us to present a solution of time complexity $O(n^{3})$ for the
    conjugacy problem for partially commutative groups  (see Theorem \ref{complexityofconprob})
    and the algorithm in Section
    \ref{ss:scemeforconprob}.

\end{itemize}

In Sections \ref{sec:DIV} and \ref{sec:conprob} we describe the
algorithms for solving all major problems in partially commutative
groups: divisibility problem, computation of the greatest common
divisor and the least common multiple of a pair of elements,
computation of maximal divisor related to fixed subgroup,
computation of block decomposition for an element of partially
commutative group, algorithms for computation of normal forms and
solving the conjugacy problem.

The principle motivation for the present paper was the desire to
understand why genetic  \cite{craven} and semi-heuristic
\cite{ESYP} algorithms happened to be so efficient for solving
conjugacy search problem and some other algorithmic problems on
special classes of partially commutative groups. See \cite{BBB}
for a brief discussion of genetic algorithms. The analysis shows
that the factors behind the good performance of these
semi-deterministic algorithms are manifestation of a good
divisibility theory. It also helps us to show, in the final
section of the paper, that our deterministic algorithms have at
most cubic time complexity.


\begin{figure}[p]

\hrulefill
\vspace{3ex}
\begin{tabular}{ccc}
  $X= \left\{x_1, \dots,x_r \right\}$ & --- & \parbox{9cm}{a finite alphabet}
  \\[2ex]
  $\GG_r$, or $\GG(X)$, or \GG & --- & \parbox[t]{9cm}{free partially commutative
  group of  rank $r$ generated by the set $X$} \\[3ex]
  $\Gplus_r$, or $\Gplus(X)$, or \Gplus & --- & \parbox[t]{9cm}{free partially commutative
  monoid of  rank $r$ generated by the set $X$} \\[3ex]
  $|w|$ & --- & \parbox[t]{9cm}{the length of the word $w$}\\[2ex]
  $[w]$ & --- & \parbox[t]{9cm}{the element of the group $\GG_r$ represented by the word $w$}\\[2ex]
 $u=v$  & --- & \parbox[t]{9cm}{equality of the elements in the group
                $\GG_r$}\\[2ex]
$u \simeq v$  & --- & \parbox[t]{9cm}{equality of words in the
free monoid $M(X\cup X^{-1})$}\\[2ex]
  $\lg(w)$ & --- & \parbox[t]{9cm}{the length of a geodesic word $w'$ such that $w =_{\GG_r} w'$}\\[2ex]
   $w=g \circ h$  & --- & \parbox[t]{9cm}{cancellation-free
  multiplication, that is, \[\lg(w)=\lg(g)+\lg(h).\]}\\[2ex]
  $\n(w)$ & --- & \parbox[t]{9cm}{the normal form of $w$}\\[2ex]
  $RC(w)$ & --- & \parbox[t]{9cm}{cyclically reduced conjugate to $w$ }\\[2ex]
  $\az(w)$ & --- & \parbox[t]{9cm}{the set
                        of symbols which occur in the word $w$}\\[2ex]
 $ \BA(w)$ & --- & \parbox[t]{9cm}{the subgroup of $\GG$ generated
 by all symbols which commute with $w$ and do not belong to $\az(w)$}\\[3ex]
 $D_l(w)$ & --- & \parbox[t]{9cm}{the set of left divisors of $w$}\\[2ex]
  $u \mid v$ & --- & \parbox[t]{9cm}{left divisibility, $u$ divides $v$}\\[2ex]
 $\ad(w)$ & --- & \parbox[t]{9cm}{the the maximal (left) abelian divisor of $w$}\\[2ex]
  $w\left[ i,j \right]$ & --- & \parbox[t]{9cm}{the interval of the word $w$ from the $i$-th
  letter and up to and including the $j$-th one, $i \leqslant j$}\\[3ex]
  $w\left[ i \right]$ & --- & \parbox[t]{9cm}{the $i$-th letter in the word
  $w$}\\[2ex]
  $e_Y(w)$ & --- & \parbox[t]{9cm}{the exhausted form of $w$
  with respect to the set $Y \subseteq X$}  \\[2ex]
  $\gd(u,v)$  & --- & \parbox[t]{9cm}{the greatest common left (right) divisor of the
  geodesic words $u$ and $v$}\\[3ex]
  $\lm(u,v)$  & --- & \parbox[t]{9cm}{the left (right) least  common multiple  of the
  geodesic words $u$ and $v$}\\[3ex]
 $p_v(w)$  & --- & \parbox[t]{9cm}{the greatest left (right)
divisor of $w$ such  that $\az(v)$ and $\az(p_v(w))$ commute
elementwise,
$\left[ \az(v), \az(p_v(w)) \right]=1$}\\[3ex]
  $\gd_Y(w)$  & --- & \parbox[t]{9cm}{the greatest  left (right) divisor of
  the geodesic word $w$  which belongs to the subgroup $\GG(Y)$}\\[2ex]

\end{tabular}
\caption{The table of notation}
\vspace{1ex}
\hrulefill
\end{figure}

\section{Free partially commutative groups} \label{sec:RAG}
\subsection{Definitions}
From now on $X=  \left\{ x_1, \ldots , x_r \right\}$ always stands
for a finite alphabet, its elements being called \emph{symbols}.
In what follows,  $x_i, x_j, \dots, y_i, y_j, \dots,$ etc.\ are
reserved for symbols from $X$ or their inverses.

Recall that a (free) \emph{partially commutative} group is a group
given by generators and relations of the form
\begin{equation} 
\GG= \left< X\mid R_X \right> \hbox{ for } R_X \hbox{ a subset of
} \left\{ \left[ x_i, x_j \right] \mid x_i, x_j \in X \right\}
\end{equation}
 This means that
the only relations imposed on the generators are commutation of
some of the generators. The number $r= |X|$ is called the
\emph{rank} of $\GG$. In particular, the free abelian group on $X$
is a partially commutative group. We also denote \[R= \left\{
\left[ x_i, x_j \right] \mid x_i, x_j \in X \right\},\] and, for
$Y \subset X$,
\[
R_Y= \left\{  \left[ x_i, x_j \right]\mid \; x_i,x_j \in Y
\right\} \cap R_X.
\]

 Set $X_{r-1} = \{\,x_2,\dots, x_r\,\}$. It is easy to see that every partially commutative  group $\GG_r=
\left< X_r\mid R_r \right>$ is an $HNN$-extension of $\GG_{r-1}=
\left< X_{r-1}\mid R_{r-1} \right>$ by the element $x_1$. In this
extension, associated subgroups $H$ and $K$ coincide and are
subgroups of $\GG_{r-1}$ generated by all symbols $x_j\in X$ such
that $x_1$ commutes with $x_j$ and $j \ne 1$. The corresponding
isomorphism $\phi: H \rightarrow K$ is the identity map ${\rm id}:
H \rightarrow H$. Indeed, this can be seen from the obvious
presentation for $\GG$:
$$
\GG= \left< \GG_{r-1}, x_1 \mid {\rel}(\GG_{r-1}), \left[x_1, x_j
\right]= 1 \hbox{ for } \left[x_1, x_j \right] \in R_X \right>
$$

Two extreme cases of this construction are worth mentioning. If
$H= K$ is the trivial subgroup then $G$ is the free product \[\GG=
\GG_{r-1}*\left< x_1 \right>\] of two partially commutative
groups. From the algorithmic point of view, free products of
groups behave nicely and most decision problems can be reduced to
similar problems for the factors.

Another extreme case, $H= K= \GG_{r-1}$, yields the direct product
\[\GG= \GG_{r-1} \times \left< x_1 \right>.\] Again,  all algorithmic
problems considered in this paper can be easily reduced to the
corresponding problems for the direct factors.

\begin{prop}\
\begin{enumerate}
\item Free and direct  products of partially commutative  groups
are also partially commutative  groups. \item Let $Y \subset X$
and $\GG(Y) = \left< Y\right>$. Then $\GG(Y)= \left< Y \mid R_Y
\right>$. Therefore the group $\GG(Y)$ is also a partially
commutative group.
\end{enumerate}
\end{prop}

The subgroups $\GG(Y)$ for $Y \subset X$ will be called
\emph{parabolic}.

\begin{proof} (1) follows directly from the definitions of free and
direct products. For the proof of (2) assume that $X=Y \sqcup
\left\{\, z_1, \dots , z_k \,\right\}$. Notice that $\GG(X)$ is an
$HNN$-extension of the group $\GG(X\smallsetminus \left\{ z_1
\right\})$, while $\GG(X\smallsetminus \left\{ z_1 \right\})$ is
an $HNN$-extension of the group $\GG(X\smallsetminus \left\{ z_1,
z_2 \right\})$, etc. Now (2) follows by induction.
\end{proof}

Notice that $\GG(Y)$ is a \emph{retract} of \GG, that is, the
image of $\GG$ under the idempotent homomorphism
\bea \GG &\rightarrow& \GG(Y)\\
y &\mapsto& y \qquad\hbox{ if } y \in Y\\
x &\mapsto& 1 \qquad\hbox{ if } x \in X \smallsetminus Y
 \eea

\subsection{\textsc{ShortLex} and normal forms}

\subsubsection{Words and group elements}
We work mostly with words in alphabet $X \cup X^{-1}$ and
distinguish between the equality of words which denote by symbol
$\simeq$, and the equality of elements in the group $\GG$, which
we denote $=$. However, we use the convention: when $w$ is a word
and we refer to it as if it is an element of $\GG$, we mean, of
course, the element $[w]$ represented by the word $w$. We take
special care to ensure that it is always clear from the context
whether we deal with words or elements of the group.

\subsubsection{\textsc{ShortLex}} Our first definition of normal forms involves the
\textsc{ShortLex} ordering. We order  the symmetrised alphabet $X
\cup X^{-1}$ by setting
\[x_1 \leqslant x_1^{-1} \leqslant x_2 \leqslant \dots \leqslant x_r \leqslant x_r^{-1}.\] This ordering gives
rise to the \textsc{ShortLex} ordering $\le$ of the set  of all
words (free monoid) in the alphabet $X \cup X^{-1}$:  we set  $u
\leqslant v$ if and only if either
\begin{itemize}
    \item  $|u| < |v|$, or
    \item $|u|=|v|$ and $u$ precedes $v$ in the sense of the
    right lexicographical order.
\end{itemize}

\textsc{ShortLex} is a total ordering of $M(X \cup X^{-1})$. If
$w$ is a word and $\left[ w \right]$ is the class of all words
equivalent to $w$ in $\GG_r$, we define $w^*$ as the $\leqslant$-minimal
representative of $\left[ w \right]$.
 The element $w^*$ is called the {\rm
(}\textsc{ShortLex}{\rm )} \emph{normal} form of the word $w$.

\subsubsection{The Bokut-Shiao rewriting rules} \label{Bokut}
 Bokut and Shiao \cite{Bokut} found a
complete set of Knuth-Bendix rewriting rules for the group
$\GG_r$:

\begin{eqnarray} \label{BS3}
  x_i^{\epsilon} x_i^{-\epsilon} &\rightarrow& 1\\
x_i^{\epsilon} w(x_{k_1},\dots, x_{k_t}) x_j^{\eta} &\rightarrow&
x_j^{\eta} x_i^{\epsilon}w(x_{k_1},\dots, x_{k_t}) \label{BS4}
\end{eqnarray}
Here, as usual, $\epsilon, \eta = \pm 1$. The first rule is
applied for all $i=1, \dots, r$,    the second one for all symbols
$x_i$, $x_j$ and $x_{k_1},\dots, x_{k_t}$ and all words
$w(x_{k_1},\dots, x_{k_t})$ such that $i
> j > k_s$ while $\left[x_i, x_j\right]=1$,
$\left[x_j, x_{k_s}\right]=1$ and $\left[x_i, x_{k_s}\right] \ne
1$.

In the case of free partially commutative monoids the rules
(\ref{BS4}) suffice, cf.\ Anisimov and Knuth \cite{AK}.

 Note that these transformation are decreasing with respect to  the \textsc{ShortLex} order.
 Bokut and Shiao showed  every
word $w$ can be transformed  into its  \textsc{ShortLex} normal
form $w^*$ using a sequence of transformations (\ref{BS3}) and
(\ref{BS4}),
\begin{equation}
\label{BS5} w=w_1 \rightarrow \dots \rightarrow w_k=w^*.
\end{equation}

\subsubsection{Geodesic words and the Cancellation Property}
A \emph{geodesic form} of a given element $w \in \GG_r$ is a word
that has minimal length among all words representing $w$. A word
is \emph{geodesic} if it is a geodesic form of the element it
represents.

\begin{lem}[Cancellation Property] \label{lem:new}
Let $w$ be a non-geodesic word. Then there exists a subword
$yw[i,j]y^{-1}$ where $y \in X \cup X^{-1}$ and $y$ commutes with
every symbol in the interval $w[i,j]$
\end{lem}
\begin{proof}
The lemma follows from  rewriting rules (\ref{BS3}) and
(\ref{BS4}) by induction on the length  $k$ of a sequence of
transformations (\ref{BS5}).
\end{proof}

\begin{lem}[Transformation Lemma] \label{lem:geodes1}
Let $w_1$ and $w_2$ be geodesic words which represent the same
element of\/ $\GG_r$. Then we can transform the word\/ $w_1$ into
$w_2$ using only the commutativity relations from $R_X$.
\end{lem}

Such transformations of geodesic words will be called
\emph{admissible transformations} or \emph{admissible
permutations}.

\begin{proof} First of all, observe that the \textsc{ShortLex} normal form
is geodesic. Indeed, the transformation of word into a normal form
does not increase the length of the word. Now, either $w_1$ or
$w_2$ can be transformed to normal form using rewriting rules
(\ref{BS3}) and (\ref{BS4}). But $w_1$ and $w_2$ are geodesic,
therefore the lengths of $w_1$ and $w_2$ can not decrease.
Since the rule (\ref{BS3}) decreases the length, this means that
we use only the commutativity rules (\ref{BS4}), which, in their
turn, can be obtained by as a reversible sequence of application
of commutativity relations from $R_X$.
\end{proof}

\subsubsection{Consequences of the Transformation Lemma}
Our first observation is that the monoid $\Gplus$ generated in
$\GG$ by the set $X$ (without taking inverses!) is the free
partially commutative monoid on the set $X$ in the sense that it
is given by the same generators and relations as $\GG$, but in the
category of semigroups.

With the help of Transformation Lemma \ref{lem:geodes1}, many
results in this paper will be proven, and definitions given, first
for geodesic words, and then transferred to elements in $\GG_r$
because they do not depend on a particular choice of a geodesic
word representing the element. For example, if $w$ is a geodesic
word, we define $\az(w)$ as the set of symbols occurring in $w$.
It immediately follows from Lemma \ref{lem:geodes1} that if $v$ is
another geodesic word and $u = w$ in $\GG_r$ then $\az(w) =
\az(v)$. This allows us to define sets $\az(w)$ for elements in
$\GG_r$. Notice that

\begin{lem} \label{lem:clem}
\label{lm:A=I} A symbol $x \in X$ commutes with $w \in \GG_r$  if
and only if $x$ commutes with every symbol $y\in \az(w)$.
\end{lem}

We set $\BA(w)$ to be the subgroup of $\GG_r$ generated by all
symbols from $X \smallsetminus \az(w)$ that commute with $w$ (or,
which is the same, with $\az(w)$).

Notice that we call elements of our alphabet $X$ \emph{symbols}.
We reserve the term \emph{letter} to denote an occurrence of a
symbol in a geodesic word. In a more formal way, a letter is a
pair (symbol, its placeholder in the word). For example, in the
word $x_1x_2x_1$ the two occurrences of $x_1$ are distinct
letters. If $Y \subset X$ and $x$ a letter, we shall abuse
notation and  write $x\in Y$ if the symbol of the letter $x$
belongs to $Y$. This convention allows us to use the full strength
of Transformation Lemma \ref{lem:geodes1} and talk about
application of rewriting rules (\ref{BS5}) to geodesic words as a
rearrangement of the letters of the word. It is especially helpful
in the discussion of divisibility (Section~\ref{sec:DIV}).

The following two results are obvious corollaries of the
Transformation Lemma.

\begin{lem} \label{lm:relative-order}
A word $u$ can be transformed to any of its geodesic forms using
only cancellations {\rm (\ref{BS3})} and permutations of letters
allowed by the commutativity relations from $R_X$, that is, rules
{\rm (\ref{BS4})} and their inverses.
\end{lem}

We shall call such transformations of words \emph{admissible}.

\begin{lem} \label{lm:cancellations}
If $x$ and $y$ are two letters in a geodesic word $w$, and
$[x,y]\ne 1$, than admissible permutations of letters in $w$ do no
change the relative position of $x$ and $y$: if $x$ precedes $y$,
than it does so after a permutation.
\end{lem}

\begin{lem} \label{lm:cancellation}
If $u$ and $v$ are geodesic words. The only cancellations which
take place in the process of an admissible transformation of the
word $uv$ into a geodesic form are those which involve a letter
from $u$ and a letter from $v$.
\end{lem}

\begin{proof} Colour letters from $u$ red and letters from $v$
black. Assume the contrary, and suppose that two black letters $x$
and $y$ cancel each other in the process of admissible
transformation of $uv$ into a geodesic form. We can chose $x$ and
$y$ so that they are the first pair of black letters cancelled in
the process. Assume also that $x$ precedes $y$ in $v$. Since the
word $v$ was geodesic, the letters $x$ and $y$ were separated in
$v$ by a letter, say $z$, which does not commute with the elements
$x$ and $y= x\inv$. Hence $z$ has to be cancelled out at some
previous step, and, since $x$ and $y$ is the first pair of black
letters to be cancelled, $z$ is cancelled out by some red letter,
say $t$. But $t$ and and $z$ are separated in the original word
$uv$ by the letter $x$ which is present in the word at all
previous steps, preventing $t$ and $z$ from cancelling each other.
This contradiction proves the lemma.
\end{proof}

\section{Divisibility Theory}
\label{sec:DIV}

In this section we transfer the concept of divisibility from
commutative algebra to the non-commutative setting of partially
commutative  groups. We follow the ideas of divisibility theory
for the
 positive Artin semigroup \cite{BrSa}. However, we shall soon see
 that, in a simpler context of partially commutative  groups, the whole
 group admits a good divisibility theory.

\subsection{Divisility} We use the notation $ v \circ w$ for multiplication in
the group $\GG_r$ when we wish to emphasise that there is no
cancellation in the product of geodesic words representing the
elements $v$ and $w$, that is, $l(v \circ w)= l(v)+l(w)$. We shall
use it also to denote concatenation of geodesic words (when the
result is a geodesic word) We skip the symbol `$\circ$' when it is
clear from context that the product is cancellation free.

We start by  introducing the notion of a \emph{divisor} of an
element; it will play a very important role in our paper.

\begin{defn} \label{defn:div}
Let $u$ and $w$ be elements in $\GG_r$. We say that $u$ is a
left\/ {\rm (}right{\rm )} \emph{divisor} of $w$ if there exists
$v\in \GG_r$ such that $w= u \circ v$ {\rm (}$w=  v\circ u$,
respectively{\rm )}.
\end{defn}

Notice that, in the context of free partially commutative monoids,
divisors have been called \emph{prefixes} \cite{BGMS}.

We set $D_l(w)$ (correspondingly $D_r(w)$) to be the set of all
left (correspondingly right) divisors of  $w$. We shall also use
the symbol $\quad\mid\quad$ to denote the left divisibility: $u
\mid v$ is the same as $u \in D_l(v)$. The terms  \emph{common
divisor}  and  \emph{common multiple} of two elements are used in
their natural meaning.

We shall work almost exclusively  with the left divisibility,
keeping in mind that the case of right divisibility is analogous.
Indeed, the two concepts are interchangeable by taking the
inverses of all elements involved.

\begin{lem} \label{lem:geodes}
Let $v$ and $w$ be fixed geodesic forms of two elements from
$\GG_r$. Then $v \mid w$ if and only if there are letters
$x_1,\dots, x_k$ in $w$ such that
 \bi
\item  $x_1\cdots x_k = u$. \item The letters $x_1,\dots, x_k$ can
be moved, by means of commutativity relations from $R_X$ only, to
the beginning of the word producing the decomposition \[ w =
x_1\cdots x_k \circ w'. \]  \item This transformation does not
change the relative position of other letters in the word, that
is, if $y_1$ and $y_2$ are letters in $w$ different from any of
$x_i$, and $y_1$ precedes $y_2$ in $w$ then $y_1$ precedes $y_2$
in $w'$.
\item
 Moreover, if we colour letters $x_1,\dots, x_k$ red and the
 remaining letters in $w$ black, then a sequence of admissible
 transformations of $w$ to $x_1\cdots x_k \circ w'$ can be chosen
 in such way that first we made all swaps of letters of different
 colour and then all swaps of letters of the same colour.
  \ei
\end{lem}

\begin{proof} This is a direct consequence of Transformation Lemma \ref{lem:geodes1}.
\end{proof}

In view of Lemma~\ref{lem:geodes}, checking divisibility of
elements amounts to simple manipulation with their geodesic forms.

\begin{lem} \label{lm:div=geodesic} The following two conditions are equivalent:
\bi \item[(a)] $u \mid w$

\item[(b)] $\lg(u\inv w) = \lg(w)-\lg(u)$. \ei
\end{lem}

\begin{proof} By definition, $u \mid w$ means that $w = u \circ v$
for some $v$, that is, $\lg(w) = \lg(u)+\lg(v)$. But then, of
course, $\lg(u\inv w) = \lg(u \inv \cdot u\cdot v) = \lg(v)$,
proving(b).

Next, assume (b). We work with some fixed geodesic words
representing $u$ and $w$. For the equality $\lg(u\inv w) =
\lg(w)-\lg(u)$ to be satisfied, $2\lg(u)$ letters need to be
cancelled in the word $uw$. By   Lemma~\ref{lm:cancellations},
every cancellation involves a letter from $u\inv$ and a letter
from $w$, which means that $u\inv$ is entirely cancelled.

Let us paint letters from $u$ red, letters from $w$ black, and
repaint green those black letters which are to be cancelled by red
letters. Obviously, if $x$ is a green letter in $w$, then the
letters in $w$ to the left of it are either green or, if black,
commute with $x$. Hence all green letters can be moved to the
left, forming a word $u'$ such that $w = u'v'$. Moreover, this
transformation can be done in such a way that the relative order
of green letters and the relative order of black letters are
preserved. Now it is obvious that $u\inv u' =1$ and $v' = u\inv
w$. Hence $u'=u$ and $w = u \circ v'$, and, in other words, $u
\mid w$.
\end{proof}

Since the rewriting rules (\ref{BS3}) and (\ref{BS4}) give an
algorithm for computing normal forms of elements,
Lemma~\ref{lm:div=geodesic} yields an obvious algorithm which
decides, for given elements $u$ and $w$, whether $u \mid w$.

\subsection{Abelian divisors and chain decomposition} We call an element $u \in \GG_r$
\emph{abelian} if all symbols in $\az(u)$ commute.

\begin{lem} \label{lem:chain}
Let $w \in \GG$. Then there exists the greatest {\rm (}left{\rm )}
abelian divisor $u$ of $w$ in the sense that if $v$ is any other
abelian left divisor of $w$ then $v \mid u$.
\end{lem}

We shall denote the greatest left abelian divisor of $w$ by
$\ad(w)$.

\begin{proof} $\ad(w)$ is composed of all letters in a geodesic form of $w$
which can be moved to the leftmost position (hence commute with
each other). Notice that  $\ad(w)$ does not depend on the choice
of a geodesic form of $w$.
\end{proof}

Lemma~\ref{lem:chain} gives us a useful decomposition of elements
in $\GG_r$. Set $w= w_1$ and decompose \bea
 w_1 & = & \ad(w_1) \circ w_2\\
 w_2 & = & \ad(w_2) \circ w_3\\
 & \vdots & \\
\eea
Iteration of this procedure gives us a decomposition of $w$:
\begin{equation} \label{eq:decomp}
w= c_1 \circ c_2\circ \cdots\circ c_t, \ \hbox{here} \ c_i =
\ad(w_i).
\end{equation}
The elements $c_i$, $i= 1, \dots , t$ are called \emph{chains} of
$w$, and decomposition (\ref{eq:decomp}) \emph{chain
decomposition}.

\begin{alg}[Chain Decomposition]
\label{alg:chain} {\rm\  \bi \item[\textsc{Input:}] a word $w$ of
length $n$ in a geodesic form. \item[\textsc{Output:}] Abelian
words $c_1,\dots, c_n$ such that
\[
w = c_1 \circ \dots \circ c_n\] is the chain decomposition of $w$.
\item[1\1] \textsc{Initialise} \bi \item[1.1.] Words $c_1 = 1,
\dots, c_n =1$ (chains). \item[1.2.] Integers $m[x] = 0$ for $x
\in X$ (chain counters). \ei

\item[2\1] \textsc{For} $i = 1,\dots, n$ \textsc{do} \bi

\item[2.1.] Read the letter $w[i]$ and its symbol $x$, $\{\, x\,\}
= \az(w[i])$.

\item[2.2.] Make the set of symbols $Y \subseteq X$ which do not
commute with $x$.

\item[2.3.] Compute the chain number for $w[i]$:
\[
m[x] = 1 + \max_{y\in Y} m[y].
\]

\item[2.4.] Append the chain:
\[
c_{m[x]} \leftarrow c_{m[x]} \circ w[i].
\]

 \ei
 \item[\textsc{End}]
 \ei
}
\end{alg}

Notice that this algorithm works in time $O(rn)$.

\begin{expl} \label{expl:nf}
{\rm
 The chain decomposition of the element $x_2 x_5
x_1x_3^{-1}x_1x_5x_4$ in the group
\[ G = \left< x_1,\dots, x_5 \mid x_ix_j = x_jx_i\; \hbox{ for }\; |i-j| >1 \right>\]
is \[x_2 x_5 x_5 \cdot x_1 x_3^{-1}x_1  \cdot x_4.\] This can seen
from the following table which captures the working of
Algorithm~\ref{alg:chain}. Here, the integers under the letters
are current values of counters $m[x]$.
\[
\begin{array}{ccccccc}
x_2 & x_5 & x_1 & x_3^{-1} & x_1 & x_5 & x_4\\
1   & 1   & 2   & 2        & 2   & 1   & 3
\end{array}
\]
}
\end{expl}

\subsection{The greatest common divisor and the least common multiple}

Obviously, the relation `$u \mid w$' is a partial ordering of
$D_l(w)$. We shall soon see that $D_l(w)$ is a lattice. (This
result is well known in the case of free partially commutative
monoids \cite[p.~150]{BGMS}.) This will allow us to use the terms
\emph{greatest common divisor} $\gd(u,v)$ and  \emph{least common
multiple} in their usual meaning when restricted to $D_l(w)$.
Without this restriction, the least common multiple of two
elements $u$ and $v$ does not necessary exist; consider, for
example, elements $u = x_1$ and $v = x_2$ in the free group $F_2 =
\left< x_1, x_2 \right>$. Notice, however, that the greatest
common divisor exists for any pair of elements in any right angled
group. To see that, concider geodesic forms of elements $u$ and
$v$; the greatest common divisor of $u$ and $v$ is formed by all
letters in $v$ which are cancelled in the product $u^{-1}v$, see
Lemma~\ref{lm:cancellation}.

The structure of the partially ordered set $D_l(w)$ in the case of
free abelian groups is quite obvious:



\begin{expl} \label{expl:2-2}
{\rm
 Let $\GG_r= \ZZ^r$ be the free abelian group of rank $r$.
Assume $w= x_1^{\alpha_1} \cdots x_r^{\alpha_r}$ written in normal
form. Then every divisor of $w$ has the form $x_1^{\beta_1} \cdots
x_r^{\beta_r}$, where $\sign(\beta _i)= \sign(\alpha _i)$ and
$|\beta _i| \leqslant |\alpha _i|$, $i= 1, \ldots , n$. In that
special case, the set $D_l(w)$ is obviously a lattice. Indeed, if
$u= x_1^{\beta_1} \cdots x_r^{\beta_r}$ and $v= x_1^{\gamma_1}
\cdots x_r^{\gamma_r}$ are elements in $D_l(w)$ represented by
their normal forms then
\bea \gd(u,v)  &=& x_1^{\phi_1} \cdots x_r^{\phi_r} \\
&& \phi_j=  \min \left\{ |\beta_j|, |\gamma_j| \right\},\quad \sign(\phi _i)= \sign(\alpha _i) \\
 \lm(u,v) &=& x_1^{\psi_1}
\cdots x_r^{\psi_r},\\
 && \psi_j = \max \left\{ |\beta_j|, |\gamma_j| \right\}, \quad
\sign(\psi _i)= \sign(\alpha _i) \eea }
\end{expl}

\begin{lem} \label{lem:nonum}
Assume that $u \mid v$ and let $u=u_1\cdots u_k$ and $v=v_1\cdots
v_l$ be chain decompositions of $u$ and $v$. Then  $k\leqslant l$
and
 $u_i\mid v_i$ for all $i\leqslant k$.
\end{lem}

\begin{proof}
Since $u_1$ is an abelian divisor of $v$ and $v_1$ is the greatest
abelian divisor of $v$, we see that $u_1 \mid v_1$.

Next decompose $u=u_1 \circ u'$ and $v=v_1 \circ v'$. Notice that
none of the letters of $u'$ appear in $v_1$, for otherwise that
would have to be included in $u_1$. Hence $u' \mid v'$ and we can
conclude the proof by induction on $\lg(u)$.
\end{proof}

We can reverse the statement of Lemma~\ref{lem:nonum} in the
important case when  $u$ and $v$ have a common multiple.

\begin{lem} \label{lem:chain-div-reverse}
Assume that $u \mid w$ and $v \mid w$. Take the chain
decompositions $u=u_1\cdots u_k$ and\/ $v=v_1\cdots v_l$ of\/ $u$
and $v$. Assume that $u_i \mid v_i$ for all $i \leqslant k$. Then
$u \mid v$.
\end{lem}

\begin{proof} Let $w = w_1\cdots w_t$ be the chain decomposition
of $w$. By Lemma~\ref{lem:nonum}, $u_1 \mid w_1$ and $v_1 \mid
w_1$. Denote $u' = u_2\cdots u_k$, $v' = v_2 \cdots v_l$ and $w' =
w_2\cdots w_t$. Arguing as in the proof of Lemma~\ref{lem:nonum},
we see that no letter of $u'$ or $v'$ appears in $w_1$ and $u'
\mid w'$ and $v' \mid w'$. We can conclude by induction that $u'
\mid v'$ and that letters of $u'$ can be chosen to lie in $v'$.
Since $u_1$, $v_1$, $w_1$ are abelian divisors, we can conclude
without loss of generality that letters of $u_1$ lies in $v_1$.
Now it is obvious that all letters of $u$ can be chosen to belong
to $v$ and that $u \mid v$.
\end{proof}

\begin{prop} \label{prop:lattice} $D_l(w)$ is a lattice with respect to the divisibility
relation $\;\mid\;$.
\end{prop}

\begin{proof} This immediately follows from Lemma~\ref{lem:nonum} and Lemma~\ref{lem:chain-div-reverse}.
\end{proof}

\begin{lem} \label{defn:gcdlcm}
Assume that $u \mid w$ and $v \mid w$. Take the chain
decompositions $u=u_1\cdots u_k$ and\/ $v=v_1\cdots v_l$ of\/ $u$
and $v$. Assume that $k \leqslant l$. If\/ $k \ < l$, set $u_{k+1}
= \dots =u_l =1$. In this notation,
\[ \gd(u,v)= \gd(u_1,v_1) \cdots \gd(u_k,v_k) \] and
\[\lm(u,v)= \lm(u_1,v_1) \cdots \lm(u_l,v_l).\]
\end{lem}

\begin{proof} This also immediately follows from Lemma~\ref{lem:nonum}
and Lemma~\ref{lem:chain-div-reverse}.
\end{proof}

\subsection{Parabolic divisors}
We set $D_{l,Y}(v)$ be the set of all left divisors of $v$ from
$\GG(Y)$; we shall call them \emph{parabolic divisors}.

\begin{prop} \label{prop:24}
For an arbitrary element $v \in \GG_r$ and subset $Y \subseteq X$
the set $D_{l,Y}(v)$ is a lattice. In particular, $D_{l,Y}(v)$
possesses unique maximal element $\md_Y (v)$.
\end{prop}
\begin{proof} This is an immediate corollary of Proposition~\ref{prop:lattice}.
 \end{proof}

\subsection{Divisibility: further properties}
We record an easy cancellation property:

\begin{lem} \label{lem:usdiv2}
Assume that $u \mid w$. If $u= u_1 \circ u_2$ and\/ $w= u_1 \circ
w_2$, then $u_2 \mid w_2$.
\end{lem}

\begin{proof} Rewrite the condition $u \mid w$ as $w = u \circ w_1$, then, obviously, the product
$w = u_1 \circ u_2 \circ w_1$ is also cancellation free. But we
are also given that $w= u_1 \circ w_2$. It follows that $u_2 \circ
w_1 = w_2$, hence  $u_2 \mid w_2$.
 \end{proof}



\begin{lem} \label{lem:usdiv}
Let $x_i, x_j \in X \cup X^{-1}$.
 If\/ $x_i \mid
w$ and $x_j \mid w$ then $x_i \ne x_j^{-1}$. If, in addition, $x_i
\neq x_j$ then $x_ix_j= x_jx_i$ and\/ $x_ix_j\mid w$.
\end{lem}

\begin{proof}  We work with a fixed geodesic form of $w$. It
easily follows from Lemma \ref{lem:geodes} that
\[
w \simeq u' \circ x_i \circ v'
\]
and
\[
w \simeq u'' \circ x_j \circ v''
\]
while $x_i$ commutes with every letter in the interval $u'$ and
$x_j$ commutes with every letter in the interval $u''$. But  one
of the letters $x_i$, $x_j$ precedes the other in $w$, therefore
$x_i$ and $x_j$ commute and they both can be moved to the
beginning of the word $w$ using only the commutativity relations.
If $x_i = x_j^{-1}$ then they can be cancelled, which contradicts
our assumption that $w$ is in geodesic form. Hence $x_ix_j$ is a
geodesic word and therefore   $x_ix_j \mid w$.
\end{proof}

Another useful lemma (an analogue of this lemma for positive Artin
semigroups can be found in \cite{BrSa}) is

\begin{lem}[Reduction Lemma] \label{lem:redlemma}
Let $v, w \in \GG_r$  and $x_i, x_j \in X \cup X^{-1}$ satisfy
$x_i\circ v= x_j\circ w$. Then  $x_i x_j= x_j x_i$ and there
exists $u \in \GG_r$ such that
\begin{equation} \label{eq:reduction}
v= x_j\circ u \ \hbox{and} \ w= x_i\circ u.
\end{equation}
\end{lem}

\begin{proof} Set $t = x_i\circ v= x_j\circ w$ then $x_i \mid t$ and $x_j \mid
t$. Then $x_ix_j = x_jx_i$ by Lemma~\ref{lem:usdiv}, and,
moreover, $x_ix_j \mid t$. Hence $t = x_ix_j \circ u$ for some
$u$. Hence $v = x_i^{-1}t = x_j\circ u$ and $w = x_j^{-1}t =
x_i\circ u$.
\end{proof}






\begin{lem} \label{lem:forprop}
Let $x$ be a letter and assume that $x \mid w$ and $v \mid w$.
If\/ $x \nmid v$ then $\left[ \az(v), x_i \right]= 1$.
\end{lem}
\begin{proof} As usual, we work with a geodesic form for $w$. Since  we have  $w =x \circ w_1$,
the word $x\circ w_1$ can be transformed into $v \circ w_2$ using
only the commutativity relations. If the letter $x$ is a part of
$v$ in the latter decomposition, we obviously have $x \mid v$.
Hence the letter $x$ is a part of $w_2$ and can be moved to the
leftmost position over all the letters in $v$, which means that
$x$ commutes with every symbol in $\az(v)$.
\end{proof}

\begin{lem} \label{lem:comdiv}
Let $p \mid w$ and $q\mid w$  and assume that  $\gd(p,q)=1$. Then
$\az(p) \cap \az(q) =  \emptyset$ and $\left[ \az(p),\az(q)
\right]= 1$.
 \end{lem}
\begin{proof} Let $p_1\cdot p_2 \cdots p_s$ and $q_1\cdot q_2 \cdots
q_t$ be the chain decompositions of  $p$ and $q$. Let $x$ and $y$
be letters such that $x \mid p_1$ and $y \mid q_1$. Since
$\gd(p,q)=1$ then $x \nmid q$ and $y \nmid p$. By
Lemma~\ref{lem:forprop}, $[x,\az(q)]=1$ and $[y,\az(p)]=1$. Notice
also that $x \ne y^{\pm 1}$. Hence $\az(x) \ne \az(y)$. Set $p = x
\circ p'$ and $q = y\circ q'$. Notice that $xy \mid w$; if $w = xy
\circ w'$ is the corresponding decomposition then $p'\mid w'$ and
$q'\mid w'$. If $g=\gd(p',q')\ne 1$, then, since every letter in
$g$ commutes with $xy$, we easily see that $g \mid p$ and $g \mid
q$, contrary to the assumption that $\gd(p,q)=1$. By induction,
$\az(p')\cap \az(q') =1$ and $[\az(p'),\az(q')]=1$, and the lemma
follows immediately.
\end{proof}

\begin{prop} \label{prop:comdiv}
Let $p \mid w$ and $q \mid w$; suppose $r= \gd(p,q)$ and set $p=
r\circ p'$, $q= r\circ q'$. Then $\az(p') \cap \az(q') =
\emptyset$, $\left[ \az(p'),\az(q') \right]= 1$ and $\lm(p,q)= r
\cdot p' \cdot q'=  r\cdot q' \cdot p'$.
\end{prop}

\begin{proof} It immediately follows from Lemma~\ref{lem:comdiv}. \end{proof}



\begin{defn} \label{defn:rdivr}
Let $p, u, v \in \GG_r$. We say that $p$ is a \emph{left divisor
of $u$ with respect to $v$} if and only if $p \mid u$ and $\left[
\az(p),\az(v) \right] =1$.
\end{defn}

\begin{prop}
There exists the greatest left divisor  of $u$ with respect to
$v$.
\end{prop}

We shall denote it $p_{v}(u)$.

\begin{proof} This is a special case of parabolic divisors, see
Proposition~\ref{prop:24}.
\end{proof}

\section{Normal forms arising from $HNN$ extensions}

In this section we shall develop a more efficient approach to
normal forms on partially commutative groups. It is based on a
presentation of $\GG_r$ as an $HNN$-extension of the group
$\GG_{r-1}$.

\subsection{$HNN$-normal
form} Let $w \in \GG_r$. We define the \emph{$HNN$-normal form}
$\n(w)$ by induction on $r$. If $r=1$ then $\GG_r$ is abelian and,
by definition, the $HNN$-normal form of an element is its
\textsc{ShortLex} minimal geodesic form. For the inductive step,
we assume that $r > 1$  and, for $l < r$, set $$X_l = X
\smallsetminus \{\, x_1,\dots, x_{l-1}\,\}.$$ Then  for every
group $\GG_l = \left< X_l \right>$ with $l<r$, the normal form is
already defined. Then
$$
\GG_r= \left< \GG_{r-1}, x_1 \mid {\rel}(\GG_{r-1}) \cup\{
\left[x_1, x_j \right]= 1 \hbox{ for } \left[x_1, x_j \right] \in
R_X\} \right>
$$
and the associated subgroup $A$ of this $HNN$-extension is
generated by all the symbols $x_i$ for $i > 1$ which commute with
$x_1$. By the inductive assumption, the $HNN$-normal form $\n(w)$
for a words $w$ which does not contain the symbol $x_1$ is already
defined.

It is a well-known fact in the theory of $HNN$-extensions
\cite{LynSh} that every element of the group $\GG_r$ can be
uniquely written in the form
\begin{equation} \label{eq:(345)}
w = s_0x_1^{\alpha_1} s_1x_1^{\alpha_2}s_2 \cdots
s_{k-1}x_1^{\alpha_k}v,
\end{equation}
where  $s_0, \ldots, s_{k-1}, v \in \GG_{r-1}$ and $s_i$ belong to
a fixed  system $S$ of words such that the corresponding elements
are left coset representatives of $A$ in $\GG_{r-1}$. The
parameter $k$ will be called the \emph{syllable length of} $w$.

We choose the system of representatives $S$ to satisfy the
following two conditions.
\begin{itemize}
    \item Each word $s\in S$ is written in the $HNN$-normal form in the group
$\GG_{r-1}$.
    \item If a letter $x_i \in A$ occurs  in $s$ then it can
not be moved to the rightmost position by means of the
commutativity relations of $\GG_{r-1}$.
\end{itemize}
Clearly, such system of representatives exists. If we assume, in
addition, that the word $v$ is written in the $HNN$-normal form,
then the word on the right hand side of (\ref{eq:(345)})  is
uniquely defined for every element $w\in \GG_r$; we shall take it
for the $HNN$-normal form of $w$.

\begin{prop}
For every $w\in \GG_r$, the \textsc{ShortLex} normal form $w^*$
and the $HNN$-normal form $\n(w)$ coincide.
\end{prop}
\begin{proof}
We use induction on $r$ and show that every word written in the
form $w^*$ is written in the form $\n(w)$. Let $w^{*}$ be written
in the form (\ref{eq:(345)}):
$$
w^*\simeq u_0x_{1}^{\beta_1} u_1x_{1}^{\beta_2}u_2 \cdots
u_{l}x_{1}^{\beta_l}u_{l+1}.
$$
From the theory of $HNN$-extensions we extract that $k=l$, that
$\alpha_i=\beta_i$ and that $s_i$ and $u_i$ are the elements of
the same cosets of $A$ in $\GG_{r-1}$.

Now consider the word $\n(w)$ and suppose that $\n(w)\not\simeq
w^*$. Then the word $\n(w)$ is not \textsc{ShortLex} minimal  in
the class of all words representing the element $w$. Due to the
Bokut-Shiao rewriting rules (Section~\ref{Bokut}) there exists an
interval $\n(w)\left[l,m\right]$ such that
$\n(w)\left[l\right]=x_i>x_j=\n(w)\left[m\right]$  and the letter
$x_j$  commutes with every letter of this interval. Since, by the
inductive assumption, $s_i^*\simeq \n(s_i)$ and $v^*\simeq \n(v)$,
the interval $\n(w)\left[l,m\right]$ can not be a subword of one of
these words. Therefore $\n(w)\left[l,m\right]$ involves $x_1$.
Notice that $x_i \ne x_1$, since $x_i>x_j$. But this implies that
$x_i$  commutes with $x_1$, and thus $x_i \in A$. This contradicts
the definition of the representatives $s_i \in S$, since $x_i$ is
involved in some representative $s_p$ and can be taken to the
rightmost position.
\end{proof}

If we restrict our considerations to the free partially
commutative monoid $\Gplus$ of positive words in $\GG$, then
$HNN$-normal form of elements in $\Gplus$ coincides with the
\emph{priority normal form} of Matiyasevich
\cite{DMM,Matiyasevich}.

\subsubsection{Cyclically reduced elements} We say that   $w \in \GG_r $
 is \emph{cyclically reduced} if and only if
$$
\lg(g^{-1}wg) \geqslant \lg(w)
$$
for every $g \in \GG_r$.

Observe that, obviously, for every element $w$ there exists a
cyclically reduced element conjugate to $w$.

One can find in the literature several slightly different
definitions of cyclically reduced words and elements adapted for
use in specific circumstances. In particular, a commonly used
definition of cyclically reduced forms of elements of
$HNN$-extension is given in \cite{LynSh}. We specialise it for the
particular  case
\[\GG_r= \left< \GG_{r-1}, x_1 \mid {\rel}(\GG_{r-1}), [x_1, A ]= 1
\right>.\]

\begin{defn}[$HNN$-cyclically reduced element]
An element \[ w=u_0x_{1}^{\alpha_1} u_1x_{1}^{\alpha_2}u_2 \cdots
u_{k}x_{1}^{\alpha_k}u_{k+1} \] is $HNN$-cyclically reduced if
\begin{itemize}
    \item $u_0=1$ and all $u_i \notin A$, $i=1, \dots, k$
\end{itemize}
and either
\begin{itemize}
    \item  $u_{k+1}=1$ and\/
    $\sign(\alpha_1)=\sign(\alpha_k)$, or
   $u_{k+1} \notin A$.
\end{itemize}
\end{defn}

\section{Conjugacy problem} \label{sec:conprob}

Wrathal \cite{Wrath} found an efficient algorithm for solving the
conjugacy problem in $\GG$ by reducing it to the conjugacy problem
in the free partially commutative monoid $\MM = M(X \cup X^{-1})$,
which, in its turn, is solved in \cite{LWZ} by a reduction to
pattern-matching questions (recall that  two elements $u,v \in
\MM$ are \emph{conjugate} if their exists $z \in \MM$ such that
$uz = zv$). In this section we shall give a \emph{direct} solution
to the conjugacy problem for partially commutative groups in terms
of a conjugacy criterion for $HNN$-extensions and the divisibility
theory developed in the preceding section. We also give solution
to the \emph{restricted conjugacy problem}: for two elements $u,v
\in \GG$ and a parabolic subgroup $\GG(Y)$, decide whether there
exists an element $z \in \GG(Y)$ such that $z^{-1}uz =v$, and, if
such elements $z$ exist, find at least one of them.

\subsection{Conjugacy Criterion for $HNN$-extensions}

We begin by formulating a well-known result, Collins' Lemma
\cite{LynSh}. It provides a conjugacy criterion for
$HNN$-extensions.

\begin{thm}[Collins' Lemma] \label{col}
Let $G= \left< H,t\mid t^{-1}A t=  B\right>$ be an $HNN$-extension
of the base group $H$ with associated subgroups $B$ and $A$. Let
 $$
    g= t^{\epsilon_1} h_1 \cdots h_{r-1}t^{\epsilon_r}h_r, \
g'= t^{\eta_1} h'_1 \cdots
    h'_{s-1}t^{\eta_s}h'_s
    $$
be normal forms of conjugate $HNN$-cyclically reduced elements
from $G$. Then either
\begin{itemize}
\item both $g$ and $g'$ lie in the base group $H$, or \item
neither of them lies in the base group, in which case $r= s$ and
    $g'$ can be obtained from $g$ by $i$-cyclically permuting it:
\[ t^{\epsilon_1} h_1 \cdots h_{r-1}t^{\epsilon_r}h_r \mapsto t^{\epsilon_i} h_i \cdots t^{\epsilon_r}h_r \cdot h_1 \cdots
t^{\epsilon_{i-1}} h_{i-1}
\]
and then  conjugating by an element $z$ from $A$, if $\epsilon_1=
1$, or from $B$, if $\epsilon_1= -1$.
\end{itemize}
\end{thm}

In the case of partially commutative groups we shall use the
following specialisation of  Collins' Lemma.

\begin{thm}[Conjugacy criterion in terms of $HNN$-extensions]
\label{thm:concrHNN} Let \[ G=  \left< H,t\mid  \left\{\,t^{-1}at=
a \mid a\in A\,\right\} \cup  {\rel}(H) \right> \]  be an
$HNN$-extension of the group $H$ with the single associated
subgroup $A$. Then each element of $G$ is conjugate to an
$HNN$-cyclically reduced element. There are two  alternatives for
a pair of conjugate cyclically reduced elements  $g$ and $g'$ is
from $G$:
\begin{itemize}
    \item If\/ $g \in H$ then $g' \in H$ and $g,g'$ are conjugated
    in base subgroup $H$.
    \item If $g \notin H$,  suppose that
    $$
    g= t^{\epsilon_1} h_1 \cdots t^{\epsilon_r}h_r, \ g'= t^{\eta_1}
h'_1 \cdots
    t^{\eta_s}h'_s
    $$
    are normal forms for $g$ and $g'$, respectively. Then $r= s$ and
    $g'$ can be obtained from $g$ by cyclically permuting it
    and then  conjugating by an element from $A$.
\end{itemize}
\end{thm}

As a corollary of Theorem \ref{thm:concrHNN} we derive the
following proposition.

\begin{prop} \label{prop:20}
Let $Y$ be a subset of $X$.  Then a pair of elements $g$ and $g'$
from $\GG(Y)$ are conjugate in $G =  \GG(X)$ if and only if they
are conjugate in $\GG(Y)$.
\end{prop}
\begin{proof} The group $\GG(X)$ can be obtained from $\GG(Y)$ as a multiple
$HNN$-extension with letters from $X \smallsetminus Y$. Therefore
the statement follows from the first clause of the Conjugacy
Criterion above.

Another proof follows from the observation that $\GG(Y)$ is a
retract of $\GG(X)$, that is, there is an idempotent morphism
$\pi: \GG(X) \rightarrow \GG(Y)$. If $h^{-1}gh = g'$ for some $h
\in \GG(X)$ then $\pi(h)^{-1}g\pi(h) = g'$ for $\pi(h) \in
\GG(Y)$.
\end{proof}

\begin{lem} \label{noid}
Let $w\in G$ and $y \in X \cup X^{-1}$. Then either

\bi \item $y^{-1}wy= w$, that is, $y$ commutes with every symbol
in $\alpha(y)$, or \item one of the four alternatives holds:
\begin{enumerate}
    \item $y^{-1}wy= w$, i. e. $y$ commutes with every letter involved
in $w$;
    \item $\left[ w, y \right] \ne 1$, $y \in D_l(w)$, $y^{-1}
\notin D_r(w)$ and\/ $\lg(y^{-1}wy)= \lg(w)$;
    \item $\left[ w, y
\right] \ne 1$, $y \notin D_l(w)$, $y^{-1} \in D_r(w)$ and\/
$\lg(y^{-1}wy)= \lg(w)$;
    \item $ \left[ w, y \right] \ne 1$, $y \in
D_l(w)$, $y^{-1} \in D_r(w)$ and\/ $\lg(y^{-1}wy)= \lg(w)-2$;
    \item $\left[ w, y \right] \ne 1$, $y \notin D_l(w)$, $y^{-1} \notin
D_r(w)$ and\/ $\lg(y^{-1}wy)= \lg(w)+2$.
\end{enumerate}
\ei
\end{lem}
\begin{proof}
The proof is obvious.
\end{proof}

Let $\hbox{CR}(w)$ denote the set of all cyclically reduced
elements conjugate to $w$.

\begin{prop} \label{prop:22}
Let $w$ be a cyclically reduced element  and take $v \in
\hbox{CR}(w)$. Then
\begin{itemize}
    \item[(a)] $\lg(v)= \lg(w)$, and
    \item[(b)] $v$ can be obtained from $w$ by a sequence of
    conjugations by elements from $X \cup X^{-1}$:
\[
w= v_0, v_1, \ldots, v_k = v,\] where each $v_i\in \hbox{CR}(w)$
and $v_i = y_i^{-1}v_{i-1}y_i$ for some $y_i \in X \cup X^{-1}$.
\end{itemize}
\end{prop}

\begin{proof} (a) is immediate from the definition of a cyclically
reduced element.

For a proof of (b), assume the contrary and take two distinct
elements $u$ and $v$ in $\CR(w)$ with the following properties:

\bi \item there is a geodesic geodesic word $y = x_1\cdots x_l$,
$x_i \in X \cup X^{-1}$, such that $y^{-1} u y = v$.

\item If we denote $y_i = x_1 \cdots x_i$ and set $u_i =
y_i^{-1}uy_i$, so that $u_{i+1} = x_{i+1}^{-1}u_ix_{i+1}$, $i =
0,1,\dots,l-1$, then not all $u_i$ belong to $\CR(w)$.

\item The word $y$ is shortest possible subject to the above
conditions.

\ei As usual, we work with geodesic forms of $u$ and $v$.  The
choice of $u$, $v$ and $Y$ implies, in particular, that if
$u_i=u_{i+1}$ then we can remove the letter $x_i$ from the word
$y$, which is impossible because of the minimal choice of $y$.
Therefore $u_i \ne u_{i+1}$ for all $i = 0,\dots, l-1$. Let us
colour letters from $y^{-1}$ and $y$ red and letters from $u$
black. Let $i$ be the first index such that $\lg(u_i) =
\lg(u_{i-1})$. By our choice of $u$ and $v$, $i >1$ and
$\lg(u_{i-1})) > \dots > \lg(u_1) > \lg(u)$, and by
Lemma~\ref{noid} the word $u_{i-1} = y_{i-1}^{-1} \circ u \circ
y_{i-1}$ is geodesic.

Now let us colour letters $x_i^{-1}$ and $x_i$ green and consider
the product $u_i = x_i^{-1} \cdot y_{i-1}^{-1} \cdot u \cdot
y_{i-1}\cdot x_i$. In the process of admissible transformation of
this product into a geodesic form, cancellations only happen with
pairs of letters of different colours
(Lemma~\ref{lm:cancellation}). Moreover, no cancellation between
red and black letters is possible, hence at least on of green
letters has to cancel with a red or black one.

By our choice of $i$, $\lg(u_i)\leqslant \lg(u_{i-1})$, therefore
some cancellations happen. If there is a cancellation within words
$x_i^{-1} y_{i-1}^{-1}$ and $y_{i-1} x_i$, then, due to the
inverse symmetry of these words, a symmetric cancellation also
takes place, and the resulting geodesic word has the form $u_i =
y'^{-1} u y'$ for some geodesic word $y'$ which is shorter than
$y_i$; this contradicts the minimal choice of $y$.

Assume now that there were a cancellation between a green letter
and red letters from the other side of the word, say, between
$x_i$ on the right and $y_{i-1}^{-1}$ on the left. But this means
that $x_i$ commutes with every black letter in $u$ and every red
letter in $y_{i-1}$. Hence $x_i$ commutes with every letter in
$u_{i-1}$ and $u_i = u_{i-1}$, a contradiction.

Hence a green letters cancels out a black letter; this means that
green and red letters commute and $u_i = y_{i-1}^{-1} x_i^{-1} u
x_i y_{i-1}$ and $\lg(x_i^{-1} u x_i) \leqslant \lg(u)$. Since $u
\in \CR(w)$, this measn that $x_i^{-1} u x_i \in \CR(w)$ and,
replacing $u$ by $x_i^{-1} u x_i$, we producing a pair of elements
which match our initial choice but have shorter conjugating
element $y$. This contradiction completes the proof.
\end{proof}

The notion of cyclically reduced element of partially commutative
group can be reformulated in terms of divisibility theory.

\begin{prop} \label{prop:21}
For an element $w \in \GG_r$, the following two conditions are
equivalent:

\bi \item[(a)] $w$ is cyclically reduced; \item[(b)]
 if $y \in X \cup X^{-1}$ is a left divisor of\/
$w$, then $y^{-1}$ is not a right divisor of $w$. \ei

\end{prop}

\begin{proof} (a) obviously implies (b), so we only need to prove the converse.
Assume that $w$ satisfies (b) but is not cyclically reduced.
Therefore there exists $v \in \GG_r$ such that $\lg(v^{-1}wv) <
\lg(w)$. Choose $v$ such that $\lg(v)$ is the minimal possible.
Let $v = y_1\cdots y_l$ be a geodesic form of $v$. Denote $w_0 =
w$ and  $$w_k = y_{k}^{-1}\cdots y_1^{-1}\cdot w \cdot y_1\cdots
y_k$$ and take the first value of $k$ such that $\lg(w_k) <
\lg(w_{k-1})$. As in the previous proof, colour the letters in
$v^{-1}$ red and letters in $w$ black. The minimal choice of $v$
implies that in the product $v^{-1}wv$ there is no cancellation of
letters of the same colour, and, since in the sequence of elements
$w_0,\dots,w_{k-1}$ the length does not decrease at any step, at
least one letter from each pair $\{y_i^{-1}, y_i\}$, $i=1,\dots,
k-1$, is not cancelled in the product $w_{k-1} =
y_{k-1}^{-1}\cdots y_{1}^{-1} \cdot w \cdot y_1\cdots y_{k-1}$.
Since the next element, $w_k = y_k^{-1}w_{k-1}y_k$, has smaller
length than $w_{k-1}$, this means that both letters $y_{k}^{-1}$
and $y_k$ cancel some black letters from $w$. To achieve this
cancellation, $y_{k}^{-1}$ and $y_k$ have to commute with all
previous red letters $y_1^{\pm 1},\dots, y_{k-1}^{\pm 1}$, and
hence with the corresponding black letters from $w$ which were
possibly cancelled in  $w_{k-1} = y_{k-1}^{-1}\cdots
y_{k-1}^{-1} \cdot w \cdot y_1\cdots y_{k-1}$. But this means that
$y_{k}^{-1}$ and $y_k$ are cancelled out in the product
$y_{k}^{-1} wy_k$, that is, $y_k$ is a left divisor of $w$ while
$y_{k}^{-1}$ is a right one.
 \end{proof}

Since the set $\hbox{CR}(w)$ is finite, Proposition \ref{prop:22}
provides a solution to the conjugacy problem for partially
commutative  groups.  However one can not state that this
algorithm runs in polynomial time until one finds a polynomial
bound for the cardinality of  the set $\hbox{CR}(w)$.

However, the situation is not that obvious, since the size of
$\hbox{CR}(w)$ can be  exponential in terms of the rank $r$ of the
group $\GG_r$. Consider, for example, the free abelian group $A$
freely  generated by $x_1,\dots, x_{r-1}$ and take for $\GG_r$ the
free product $A \star \left< x_r \right>$. Take $w = x_1\cdots
x_{r-1}x_r$. If $I \sqcup J = \{1,\dots,r-1\}$ is any partition of
the set $1,\dots,r-1$ into disjoint subsets $I = \{i_1,\dots,
i_k\}$ and $J= \{j_1,\dots, j_l\}$, then any element
$x_{i_1}\cdots x_{i_k} \cdot x_r \cdot x_{j_1}\cdots x_{j_l}$
belongs to $CR(w)$. This gives us $2^{r-1}$ elements in $CR(w)$.

However, the conjugacy problem in free products as one just
mentioned is trivial, which shows that the size of $CR(w)$ is not
necessarily a good indicator of its difficulty. In Section
\ref{sec:com} we give a polynomial time algorithm for solving the
conjugacy problem in partially commutative groups.

\subsection{Block decomposition}

For a partially commutative  group $\GG_r$ consider its graph
$\Gamma$. The vertex set $V$ of $\Gamma$ is a set of generators
$X$ of $\GG_r$. There is an edge connecting $x_i$ and $x_j$ if and
only if $\left[x_i, x_j \right] \ne 1$. The graph $\Gamma$ is a
union of its connected components. Assume $I_1, \ldots , I_k$ are
sets of letters corresponding to the connected components of
$\Gamma$. Then $\GG_r= \GG(I_1) \times \cdots \times \GG(I_k)$
and the words that depend on letters from distinct components
commute.

Consider $w \in \GG$ and the set $\az(w)$. For this set, just as
above, consider the graph $\Gamma (\az(w))$ (it will be a subgraph
of $\Gamma$). For this graph can be either connected or not. If it
is not, then we can split $w$ into the product of commuting words
$\left\{ w^{(j)}| j \in J \right\}$, where $|J|$ is the number of
connected components of $\Gamma(\az(w))$ and the word $w^{(j)}$
involves the letters from $j$-th connected compound. Clearly,
$\left[ w^{(j_1)},w^{(j_2)} \right] = 1$.

We shall call an element  $w\in \GG$  a \emph{block} if and only
if the graph $\Gamma(w)$ is connected.

As we already know  every element $w$ admits a block
decomposition:
\begin{equation} \label{eq:bl}
w= w^{(1)} \cdot w^{(2)} \cdots w^{(t)},
\end{equation}
 here $w^{(i)}$, $i= 1, \dots ,t$ are blocks.

\begin{prop} \label{prop:23}\
\begin{enumerate}
    \item The element $w$ of a partially commutative  group is cyclically reduced
    if and only if each block involved in decomposition
    {\rm (\ref{eq:bl})} is cyclically reduced.
    \item Let $w= w^{(1)} \cdot w^{(2)} \cdots w^{(t)}$ and
$v= v^{(1)} \cdot v^{(2)} \cdots
    v^{(s)}$ be cyclically reduced elements decomposed into the
    product of blocks. Then $v$ and $w$ are conjugate if and only
    if $s= t$ and, after some certain index re-enumeration,  $v^{(i)}$ is
    conjugate with $w^{(i)},\; i= 1, \dots ,t$.
\end{enumerate}
\end{prop}
\begin{proof} The first statement is clear. To prove the second
statement, assume that the blocks of $v$ and $w$ are conjugate.
Applying Proposition \ref{prop:20} we obtain that the blocks are
conjugated in correspondent connected components of $\Gamma (w)$
and $\Gamma(v)$. But this implies that $v$ and $w$ are conjugate.

Conversely, suppose that $v$ and $w$ are conjugate. According to
Proposition \ref{prop:22} there exists a sequence of cyclically
reduced elements $w= v_0, \ldots, v_k =  v$. Here $v_i \in
\hbox{CR}(w)$ and $v_{i+1}, v_i$ are conjugated by a letter from
$X \cup X^{-1}$, and moreover $l(v_{i+1})= l(v_i)$. Therefore the
second statement follows under induction on the length of this
sequence.
 \end{proof}

\medskip

The proposition above allows  us to reduce the proof of most
statements for cyclically reduced words in partially commutative
groups to the case when the considered word is a block.

\medskip

\subsection{Conjugacy with respect to a subgroup} \label{ss:conrelsub}

Let $\GG_r= \GG(X)$ be a partially commutative group and let $Y$
be a proper subset of $X$. We shall now draw our attention to the
question
\begin{quote}
`when are two elements of $\GG(X)$ conjugate by an element from
$\GG(Y)$?'.
\end{quote}

Assume that $w$ is a cyclically reduced block element such that
$\az(w) = X$. By $w^{\GG(Y)}$ we shall denote the subset of
elements from $\GG(X)$ conjugate with $w$ by an element from
$\GG(Y)$.

Our nearest goal is to provide an algorithm for solving the
problem of `being an element of $w^{\GG(Y)}$'.

\begin{prop} \label{prop:25}
Let be $w$ a cyclically reduced block element such that $\az(w) =
X$. Then there exists a unique element $w_0$ of $w^{\GG(Y)}$ such
that:
\begin{enumerate}
    \item $w_0$ is cyclically reduced;
    \item $D_{l,Y}(w_0)= \left\{ 1 \right\}$.
\end{enumerate}
\end{prop}
\begin{proof} \textsc{Existence.} Let $r= \md_Y(w)$, $w_1=  r \cdot w'$ and
set $w = r^{-1}wr= w' \cdot w_1$. Clearly, $D_{l,Y}(w')= \left\{ 1
\right\}$. If $D_{l,Y}(w' \cdot r)= \left\{ 1 \right\}$ then set
$w_1= w_0$.

Otherwise $r=p_1r_1$ and $p_1=p_{w'} (r)\ne 1$, $\left[ \az(p_1),
\az(w') \right]=1$. Since $\az(w') \cup Y= X$ and since $w$ is a
block element, then $\az(p_1) \varsubsetneq Y$. Set
$w_2=p_1^{-1}w_1 p_1=w'r_1p_1$.  If $D_{l,Y}(w_2)=1$ then the
process terminates. Otherwise proceed to the next iteration. Then
$p_2=p_{w'r_1}(p_1) \ne 1$ and $\az(p_2)$ is a proper subset of
$\az(p_1)$. This means that the process terminates in no more than
$|Y|$ steps. Notice that the number of steps in this procedure
depends on the cardinality of the alphabet $X$ and does not depend
on the length of the word $w$.

After iterating this procedure no more than $l(w)-1$ times we
obtain the required element.

\textsc{Uniqueness.} Suppose $w$ and $v$ are two elements of
$w^{\GG(Y)}$ satisfying conditions (1) and (2) of Proposition
\ref{prop:25}. Choose $u \in \GG(Y)$ such that $v = u^{-1}wu$. If
$v \ne w$ then $u \ne 1$. Consider geodesic forms of $u^{-1}$, $w$
and $u$ and colour their letters red, yellow and green,
correspondingly. A geodesic form for $v$ is the result of
admissible permutations and cancellations in the product $u^{-1}w
u$. If some red letter is not cancelled out in $v$ then
$D_{l,Y}(v) \ne \{1\}$. Since cancellations can take place only
between letters of different colours
(Lemma~\ref{lm:cancellation}), this means that each red letter has
to be cancelled out either by a yellow letter (but then
$D_{l,Y}(w) \ne \{1\}$, or by a green letter. In the latter case
red and green letters cancel each other, leaving the yellow letter
intact. Since the original relative order of yellow letters can be
restored by virtue of Lemma~\ref{lm:relative-order}, this means
that $v=w$.
 \end{proof}

\begin{defn} \label{defn:exh}
Assume that an element $w$ satisfies the conditions of Proposition
\ref{prop:25}. Then the unique element $w_0$ constructed in the
proposition is called the \emph{exhausted form} of $w$  with
respect to the parabolic subgroup $\GG(Y)$ and is denoted $e_Y(w)=
w_0$. We omit the subscript $Y$ when it is clear from the context.
 \end{defn}

\begin{rem} \label{rem:alg}
The proof of Proposition \ref{prop:25} provides us with an
algorithm for computation  of $e(w)$.
\end{rem}
\begin{alg}(Computation of $e(w)$) \
\begin{itemize}
    \item[1\1] Compute $w_1= \md_Y(w)$.
    \item[2\1] Find the geodesic decomposition $w= w_1w_2$.
    \item[3\1] Compute the greatest left divisor of $w_1$ with
respect to
    $w_2$: $w_3= p_{w_1}(w_2)$
    \item[4\1] Find the geodesic decomposition $w_1= w_3 \cdot w_4$.
    \item[5\1] Compute $w_5= p_{w_2 \cdot w_4}(w_3)$
    \item[6\1] Find the geodesic decomposition $w_3= w_5 \cdot w_6$.
   \item[\vdots ] \qquad\dots
    \item[\textsc{Output:}] The geodesic decomposition for $e(w)$ is
    \[e(w)= w_2 \cdot w_{4} \cdots w_{2k} \cdot w_{2k-1}, \quad 2k \leqslant     |Y|.\]
\end{itemize}
\end{alg}

\begin{defn} Let $v, w \in \GG$. We say that $v$ is a
\emph{factor} of $w$ if there exist $u_1,u_2 \in \GG$ such that
\[
w = u_1\circ v \circ u_2.
\]
\end{defn}

Notice that if $W$ is a cyclically reduced block element and
$e(w)$ is its exhausted form then it follows from the algorithm
for computation of $e(w)$ that $e(w)$ is a factor of $w^r$, where
$r = |Y|$.

\subsection{An algorithm for solving the Conjugacy Problem}
\label{ss:scemeforconprob}

 The conjugacy criterion for $HNN$-extensions (Theorem
\ref{thm:concrHNN}), Propositions \ref{prop:22}, \ref{prop:25} and
Remark \ref{rem:alg} provide an algorithm for deciding whether the
elements $g,h \in \GG_r$ are conjugate or not.

The description of procedures which are used in the steps of the
scheme is given in Section \ref{sec:com}.

\begin{alg}[Conjugacy Problem] \

\begin{itemize}
    \item[\textsc{Input}:] Two elements $g$ and $h$.
    \item[\textsc{Output}:] `Yes' if $g$ and $h$ are conjugate.
    \item[1\1] Compute the
    normal forms $g^*$ and $h^*$ for elements $g$ and $h$
    correspondingly.

    \item[2\1] Using the division algorithm and Proposition
    \ref{prop:21} compute cyclically reduced elements $g'$ and $h'$ of
    $\GG_r$ which are conjugate to $g^*$ and $h^*$ respectively.

    \item[3\1] If one of either $g'$ or $h'$ is an element of
    $\GG_{n-1}$, while the other is not, then $g$ and $h$ are not
    conjugate. If both $g'$ and $h'$ lie in $\GG_{n-1}$,
    then solve the conjugacy problem in a partially commutative  group
    $\GG_{n-1}$ of lower rank.

    \item[4\1] Both $g'$ and $h'$ are not the elements of
    $\GG_{n-1}$.
     Cyclically permuting a cyclically reduced element
    we may assume that its first letter is $x_1^{\epsilon}$, i.e. we
    obtain cyclically reduced and $HNN$-cyclically reduced word
    simultaneously.
    Regarding $g'$ and $h'$ as words written in form
    (\ref{eq:(345)}) compare their syllable lengths  $k_1$ and $k_2$. If $k_1
\neq
    k_2$, then the elements $g$ and $h$ are not conjugate.

    \item[5\1] If $k_1= k_2$, then compute the block decompositions
    $w^{(1)} \cdots w^{(t)}$ and $v^{(1)} \cdots v^{(p)}$ for $g'$ and
$h'$.
    If the number of factors in the block
    decompositions for $g'$ and $h'$ are distinct (that is, $t \neq p$),
    then, due to Proposition \ref{prop:23}, the  elements are not
conjugate.

    \item[6\1] Construct the graphs $\Gamma (g')$ and $ \Gamma (h')$.
    If the connected components of these graphs do not coincide,
     then,
    according  to Proposition \ref{prop:23}, some of the blocks $w^{(i)}$
and
$v^{(i)}$
    for  $g'$ and $h'$  are not conjugate and so are $g$ and $h$.

    \item[7\1] If $t= p>1$ and $\az(w^{(i)})= \az(v^{(i)})$ for some
block
decompositions of
    $g'$ and $h'$, then the conjugacy problem
    for $g$ and $h$ is reduced to the conjugacy problem for pairs
    of words $\left\{ w^{(i)}, v^{(i)} \right\}, i= 1, \ldots, t$ in
the     groups $\GG(\az(w^{(i)}))$ of lower ranks.

    \item[8\1]  If $t= p= 1$ and $\az(g') \subset X$, $\az(g') \neq
    X$ then the  conjugacy problem
    for $g$ and $h$ is reduced to the one in group $\GG(\az(g'))$
    of lower rank.

\item[9\1] If $t= p= 1$ and $\az(g')= X$. Let $w_i$ denote
    the $i$-cyclical permutation of the word
    $$w= s_0x_1^{\alpha_1}
s_1x_1^{\alpha_2}s_2 \cdots
    s_{k_2-1}x_1^{\alpha_{k_2}}v.$$
    According to the conjugacy criterion, it  suffices to
    solve the conjugacy problem for pairs of elements $\left\{ w_i, h'
\right\}$, $i= 0, \dots,
    k_2$ related to associated subgroup $A$.
    To complete Step 9\1 compute the
    elements $\left\{ e(w^{(i)}), e(h') \right\}$, $i= 0, \dots,k_2
    $
    (see the scheme for computation of an exhausted element in the end
of the
Section
    \ref{ss:conrelsub}).
    If every pair $\left\{ e(w^{(i)}), e(h') \right\}$ is a pair of
different
words, then
    $g'$ and $h'$ are not conjugate. If there exists an index
    $i_0$ such, that $ e(w^{(i_0)})=  e(h')$, then the words $g'$ and
$h'$ are
    conjugate in $\GG_r$ and so are $g$ and $h$.
     \item[\textsc{End}]
\end{itemize}
\end{alg}

\begin{rem}
Let $w$ be a geodesic block element. It follows from algorithm for
conjugacy problem that there might be many exhausted elements
conjugate to $w$. We consider the ShortLex-minimal among them and
call it  the \emph{canonical representative} for the class
$\left\{ w^{\GG_r}\right\}$
\end{rem}

\subsection{Double cosets of parabolic subgroups.}

Consider a partially commutative  group $G=\left< X\mid R_X
\right>$ and two parabolic subgroups $H_1=\GG(Y)$ and
$H_2=\GG(Z)$, $Y,Z \subset X$.

In this section, we outline an algorithm which decides, for given
elements $g_1$ and $g_2$ of $G$, whether they belong to the same
double coset $H_1gH_2$ or not. And if they do, we shall show how
one can easily find the elements $h_1$ from $H_1$ and $h_2$ from
$H_2$ such that $g_1=h_1g_2h_2$.

\begin{defn}
The element $g' \in H_1gH_2$ is called a \emph{canonical
representative} of the double coset $H_1gH_2$ if and only if\/
$\md_{l,Y}(g')=\md_{r,Z}(g')=1$.
\end{defn}

Using induction on the length of an element (for the elements of
$H_1gH_2$) we show that in each coset $H_1gH_2$ there exists a
canonical representative. If the element $g$ is canonical then the
algorithm terminates. Otherwise, say, $1\ne r=\md_{l,Y}(g)$ and
$g=r \circ g_1$. Then $g_1 \in H_1gH_2$ and $l(g_1)<l(g)$.
Therefore the computation of the canonical representative for
$H_1gH_2$ is equivalent to ``exhausting'' the element $g$ on the
left and on the right (see Definition \ref{defn:exh}). We
estimate the complexity of the algorithm for computation of
$\md_{l,Y}(g)$ in Section \ref{ss:complofmaxdiv}.

\begin{prop}
Every double coset $H_1vH_2$ contains a unique canonical
representative.
\end{prop}

\begin{proof}
Assume that $v$ and $w$ are two canonical representative of
$H_1vH_2$. We may assume without loss of generality that
$v\ne 1$
and $w \ne 1$; we also assume that  $v$ and $w$ are  written in geodesic form. Then $H_1vH_2 \ne
H_1H_2$. Assume that the length of $v$ is lower or equal than the
length of $w$ . Then $w=h_1vh_2$, where $h_1$ and $h_2$ are
written in geodesic form and $l(h_1)+l(h_2)=s$ is minimal.
By our choice of $h_1$ and $h_2$, no letter of $h_1$
(correspondingly, $h_2$) cancels with a letter of $h_2$
(correspondingly, $h_1$). Therefore at least one letter $y$ from
either  $h_1$ or $h_2$ cancel a letter $y^{-1}$
 in the word $v$ (see Section \ref{ss:con}). This leads to a contradiction
with the definition of a canonical representative.
\end{proof}

\medskip

\section{Complexity of algorithms: some estimates} \label{sec:com}

\subsection{Turing machines and the definition of complexity of
algorithms}

Let $G$ be a finitely presented group,
\begin{equation} \label{eq:G}
G=  \left<\, X\mid R\, \right>.
\end{equation}
The algorithms we construct in this paper are term rewriting
systems. In particular, an algorithm $\mathcal{A}$ takes a given
word $w$  in alphabet $X \cup X^{-1}$ as input and computes
another word $v= \mathcal{A}(w)$. We adopt a version of Church's
Thesis and assume that all our procedures are implemented on a
multi-taped Turing machine $M$. We denote by $t_M(w)$ the run time
of $M$ computing an output word for a fixed initial input word $w$.

As usually, we introduce the \emph{worst case complexity} of
$\mathcal{A}$ as an integer valued function
\begin{equation} \label{eq:t_M}
t_M(n)= \max \left\{\, t_M(w) \mid |w|= n\, \right\}.
\end{equation}
Here, $|w|$ denotes the length of a word $w$ in the alphabet
$X\cup X^{-1}$.

An algorithm $\mathcal{A}$ is called \emph{polynomial} if there
exists its realization on Turing machine $M$ such that $t_M(n)$ is
bounded by a polynomial function of $n$.

The second definition of complexity is preferable in  terms of
practical, experimental use. We shall define \emph{`average case
complexity'} by the following timing function:
\begin{equation} \label{eq:ac}
t_{M,C}(n)= \frac{\sum\limits_{|w|= n}t_M(w)}{|S_n|},
\end{equation}
here $|S_n|$ stands for the number of all words of length $n$.

The algorithm $\mathcal{A}$ is called \emph{polynomial on average}
if and only if there exists its realization on Turing machine $M$
such that $t_{M,C}(n)$ is bounded by a polynomial function of $n$.

It is natural and convenient to divide all algorithms
into two classes, polynomial and
non-polynomial. However, it is not
sufficiently precise when we look at the practical applications of
algorithms.

Consider, for example, a simple question. Given a letter $x$ from
an alphabet $X$ and a word $W$ in $X$, does $w$ contain an
occurrence of $x$? When the inputs are represented as sequences of
symbols on the tape of a Turing machine, this simple problem has
linear time complexity in terms of the length of $w$. Even more
counterintuitive is the complexity of the deletion of a letter
from a word: we remove the letter  and then have to close the gap,
which involves moving a segment of the word one position to the
left, rewriting it letter after letter. However, if the word is
presented as collection of triples
\begin{quote}
(letter, pointer to the previous letter, pointer to the next
letter),
\end{quote}

\noindent then deletion of a letter requires a modification of at
most three triples and can be done in constant time.

This simple example shows that assertions that certain algorithms
works in linear, quadratic, etc., time, are relative and depend on
the data formats used.

This is one of the reasons why computer scientists use a much more
pragmatic setup for the complexity of algorithms, using more natural
data formats and models of computation (for example, random access
machines) \cite{C}. We use \cite{C} as the main source of
references to standard results about complexity of algorithms used
in this paper, and, in particular, to various modifications and
generalizations of sorting algorithms. The latter are
indispensable building blocks for our word processing algorithms.
The only modification we introduce  is that in \cite{C} the
relations between members of sequences of letters (numbers) remain
constant, while in our case they are changing in the process of
work of our higher level word processing algorithms. To emphasize
the difference, we shall refer to our algorithms as
\emph{dynamic} sorting algorithms.

\subsection{Conventions} \label{ss:con}

We want to make our results independent of varying and diverse
definitions of complexity of algorithms. To that end we choose the
following tactic.

First of all, we select certain operations on
words which we regard as elementary. With a single exception, all
these operations are well known in the theory of algorithms and we
use classical results about their complexity.

Secondly, we understand complexity of an algorithm $\mathcal{A}$
 as the (worst possible) number of elementary operations
it is using when working on inputs of length $n$ and denote it
$f_{\mathcal{A}}(n)$. We formulate our results in the standard
$O$-notation: $O(f_{\mathcal{A}}(n)) =  g(n)$  if there exists
positive constant $c$ and integer $n_0$ such that
\[
0 \leqslant f_{\mathcal{A}}(n) \leqslant cg(n)
\]
for all $n > n_0$.

To introduce elementary operations, we need some notation.

Let $w$ be a word in alphabet $X \cup X^{-1}$. Denote by $l(w)$
the standard length of $w$, that is, the number of letters it
contains.  We denote as $w[i]$ the letter in the $i$-th position
of the word $w$, while $w[i,j]$ denotes the subword
$w[i]w[i+1]\cdots w[j]$.

Now we list our elementary operations. First three of them will be
also called \emph{basic} operation.

\subsubsection{Cancellation.} \label{oper1} \label{sec:basic1}
If $w[j] =  w^{-1}[j+1]$, then, after deleting letters $w[j]$ and
$w[j+1]$ from from $w$, we have a word $w'$ with $l(w') =  l(w)
-2$. The complexity of this operation is a constant \cite{C}.

\subsubsection{Transposition of a letter and an admissible
interval.}\label{oper2} \label{sec:basic2}  We shall use two types
of \emph{admissible} intervals in a geodesic word $w$. If a letter
$x = w[j]$ commutes with all letters in the interval $w[j+1,i]$,
we say that the interval $w[j+1,i]$ \emph{admits} swapping with
$x$. In that case, the result of the transposition applied to the
word
\[ \cdots w[j]w[j+1,i]\cdots\] is the word \[ w' = \cdots
w[j+1,i]w[j] \cdots .
\]

Another type of admissible intervals is the following: if $w$ is a
geodesic word and $w = xw'$, then $w'$ is an admissible interval. In
that case we can swap $w'$ and $x$.

If one uses pointers, the transposition of a letter and an
admissible interval requires only linear time (and constant time
if the length of the interval is bounded from above by some
constant $K$).

We described the left hand side version of the operation. Of
course, the swap can be made the other way round.

\subsubsection{Query} \label{oper3} \label{sec:basic3} Let $x$ be a fixed letter of our alphabet
$X$ and $w=y_1\cdots y_l$ a word, $y_i \in X \cup X^{-1}$. Then we
set
\[
\ln(w) = \left\{\begin{array}{rl}y_{i_0} & \quad\hbox{ if }\quad
[y_{i_0},x]=1 \hbox{ and } [y_i,x]\ne 1 \hbox{ for all } i < j_0,\\
1 & \quad\hbox{ if such } i_0 \hbox{ does not exist}.\end{array}
\right.
\]
The right hand side version of this operation is denoted $\rn(w)$.

\subsubsection{ Abelian sorting.} Assume that all letters in
the interval $w[i,j]$ commute pairwise (we call such intervals
\emph{abelian}). We apply to the interval one of the sorting
algorithms of \cite[Chapter~2]{C} and rewrite it in the increasing
order (with respect to indices).

Most algorithms in \cite[Chapter~2]{C} have worst case run time
$O(n^2)$ and average time complexity $O(n\ln n)$. The latter is
also a lower bound for comparison algorithms
 \cite[Section~9.1]{C}. We note that there is a
sorting algorithm working in linear time \cite[Section~9]{C}.

\subsubsection{Computation of the maximal left divisor which commutes with the letter $x$.} \label{oper5}
Given a word $w$ and letter $x$, we want to find $p =  p_x(w)$,
the maximal left divisor of $w$ which commutes with $x$, and
rewrite the word $w$ in the form $pw'$ without changing the length
of the word. We shall show in Section~\ref{sec58}, that this
operation can be performed in linear time.

Again, we can similarly introduce the right hand side version of
this operation.

\subsubsection{Iterated divisors.} Let $Z =
\{\,z_1,\dots,z_k\,\}$ be a linearly ordered (in a way possibly
different from ordering by indexes) subset of the alphabet $X$.
Define by induction
\[
p_{\{z_1,\dots,z_i\}}(w) =  p_{z_i}(p_{\{z_1,\dots,z_{i-1}\}}(w))
\]
and
\[
p_Z(w) =  p_{\{z_1,\dots,z_k\}}(w).
\]
This is an iterated version of the previous operation: finding the
maximal left divisor commuting with a given letter.

\subsubsection{The greatest common divisor and the least common
multiple of two geodesic words in the free abelian group} The
corresponding algorithm is relatively simple and clearly explained
in Example~\ref{expl:2-2}. It can be performed in linear
time.

\subsubsection{Divisors in abelian group}\label{oper8} Let $Y \subseteq X$. We
need to compute $\md_Y(w)$ in the special case when $\GG_r$ is an
abelian group. This can be done in linear time.

\subsubsection{Cyclic permutation of a geodesic word} This
operation transforms a geodesic word $w = p\cdot q$ into $x' =
q\cdot p$.

Complexity of this operation depends on the model of computation.
It is linear in terms of the basic operation~\ref{sec:basic2}.
However, if we represent the word as a sequence of pointers, then
the complexity is constant.

\subsubsection{Occurrences of letters} \label{oper10} Given a word $w$, $\az(w)$
denotes the set of letters occurring in $w$. If pointers are used,
then the set $\az(w)$ and the intersection of two such sets
$\az(u)\cap \az(v)$ can be found in constant time.

\subsubsection{Connected components of the graph $\Gamma(w)$}
\label{oper11}  It is well known that connected components of a
graph with $v$ vertices and $e$ edges can be found in linear time
$O(v+e)$ \cite[\S 23]{C}. In our context $v=r$ and $e \leqslant
r^2$, which gives a quadratic bound in terms of $r$.

\subsection{Complexity of the algorithm for computation of normal
forms} \label{sec:normal} Let $\GG_r$ be a partially commutative
group of rank $r$. If $r=1$ then $\GG_r$ is the infinite cyclic
group generated by element $x_1$. In order to reduce to a normal
form a word $w$ of length $n$ in alphabet $\{\,x_1, x_1^{-1}\,\}$,
one needs at most $[n/2]$ cancellations.

Let now $r >1$. Our alphabet $X = \{\,x_1,\dots, x_r\,\}$ is
linearly ordered as \[x_1 < \dots < x_r.\] Then $\GG_r$ is the
$HNN$-extension of the group $\GG_{r-1} = \left< x_2,\dots,
x_r\right>$ by the stable letter $x_1$ with the associated
subgroup
\[
A = \left< x_i \mid [x_i,x_1] =1, \, i \ge 2 \right>.
\]
The group $\GG_r$ can be written in generators and relations as
\[
\GG_r = \left< \GG_{r-1},x_1 \mid \rel \GG_{r-1}, \, x_1^{-1}ax_1
= a \hbox{ for all } a \in A \right>.
\]
We shall use a standard algorithm for computing normal forms of
elements of $HNN$-extensions, adapted to our concrete situation.
We shall denote it $\mathcal{N}$. Given a word $w$ of length $n$
in alphabet $X \cup X^{-1}$, it is converted to normal form
${{N}}(w)$. Let $f_r(n)$ be the number of elementary operations
needed for computation of a normal form of arbitrary word of
length $n$, and $g_r(n)$ the same quantity, but restricted to
geodesic words of length $n$.

\begin{prop}
In this notation,
$f_r(n) \leqslant n^2$,
while
$ g_r(n) \leqslant 2n.$
\end{prop}

\begin{proof}
If $r=1$, then these estimates easily follow from the estimates
for infinite cyclic group given at the beginning of this section.

Let now $r> 1$ and assume, by way of induction on $n$, that for
all $i=1,2,\dots,r-1$ the bounds
\[
f_i(n) \leqslant n^2 \hbox{ and } g_i(n) \leqslant 2n
\]
are true. Consider first the function $g_r$.  By
Lemma~\ref{lem:geodes1}, the  transformation of a geodesic word to
normal form does not change its length, that is, cancellations are
not used.

Let $w = w' \cdot x$ with $x \in X \cup X^{-1}$. Then $\n(w) =
\n(\n(w')\cdot x)$. By the inductive assumption, the transformation
of $w'$ into $\n(w')$ requires at most $2(n-1)$ elementary
operations. To prove the desired bound for $g_r$, we have to show
that transformation of $\n(w')\cdot x$ into $\n(w)$ requires at most
two elementary transformation, which can be show by easy
case-by-case considerations.

The proof of a bound for $f_r(n)$ can be done similarly. the only
difference is that we encounter the so-call \emph{pinch}, that is,
the situation when $\n(w') = \cdots vx_1^\epsilon u$, $u \in A$, $x
= x_1^{-1}$ and $v \in \GG_{r-1}$. In that case $\n(w')x = \cdots
vx_1^\epsilon u x_1^{-\epsilon}$ and we need to apply the
following elementary operations:
\begin{itemize}
\item we have to swap $u$ and $x_1^{-\epsilon}$,
\item cancel $x_1^\epsilon x_1^{-\epsilon}$,
\item apply some number of elementary operations for converting $vu$ to
normal form.
\end{itemize}
The word $vu$ is geodesic and does not belong to $A$, for
otherwise $\n(w')$ is not a normal form. Indeed, if the word $vu$
is not geodesic, then, according to
the Cancellation Property (Lemma \ref{lem:new}),
there exists an interval
$x_kvu[i+1,j]x_k^{-1}$ and letter $x_k$ which commutes with all
letters in the interval $vu[i+1,j]$, and, by our assumption, with
the letter $x_1$. Applying again Lemma \ref{lem:new},
we conclude that the word $\n(w')$ is not geodesic, a
contradiction. We see now that need at most $2(n-2)$ elementary
operations for converting $vu$ to normal form. The total number of
operations now is bounded by
\[
(n-1)^2 +2n -4 < n^2,
\]
which completes the proof.
\end{proof}

\subsection{Complexity of the algorithm for computation of the chain
decomposition} \label{sec:complexity-E'}

Assume that we have a geodesic word $w$. Recall (see
Section~\ref{sec:DIV}) that
\[
w = w_{1}\cdots w_{l}
\]
is the product of chains and $w_1$ is the largest abelian divisor
of $w$, and if we write $w = w_1\circ w'$ then $w_2$ is the
largest abelian divisor of $w'$, etc.

The largest abelian divisor of a geodesic word $w$ is computed by
the following procedure. Let $z_1$ be the first letter of $w$ and
$p_1 = p_{z_1}(w)$ the result of application of the elementary
operation of type (iv). Set $w = p_1w_1$. Applying a sorting
algorithm to the word $p_1$, we can rewrite $p_1 =
z_1^{\alpha_1}\circ p'_2$, where letter $z_1$ does not occur in
$p'_2$. Let $z_2$ be the first letter of $p'_2$. Then $p_2 =
z_1^{\alpha_1}\cdot p_{z_2}(p'_2) = z_1^{\alpha_1}
z_2^{\alpha_2}p'_3$. This process continues no more than $r$ times
and requires at most $3r$ elementary operations. The procedure
produces an abelian divisor $w_1 = z_1^{\alpha_1}\cdots
z_k^{\alpha_k}$ of and a geodesic decomposition $w = w_1\cdot w'$.
It is easy to see that $w_1$ is the largest left abelian divisor
of $w$ in the sense that any other left abelian divisor divides
$w_1$.

Denote by $h_r(n)$ the number of elementary operations needed (in
the worst case scenario) for computation of the chain
decomposition  of a geodesic word of length $n$.

\begin{prop} \qquad $h_r(n) \leqslant 3rn. $
\end{prop}

\begin{proof}
The proof immediately follows from the description of the
procedure.
\end{proof}

\subsection{Complexity of algorithms for computation of divisors}
\label{ss:complofmaxdiv} Let $u$ and $v$ be left divisors of an
element $w$. We assume that $u$ and $v$ are written in geodesic
form. We compute, with the help of the algorithm for computation
of a the chain decomposition, the first chains $u_1$ and $v_1$ in
the chain decompositions of elements $u$ and $v$:
\[
u = u_1\cdot u', \quad v = v_1\cdot v'.
\]
We compute next the set
\[
\az( u_1) \cap \az(v_1) = \{\,z_1,\dots, z_k\,\}.
\]
If this set is empty then $\md(u,v)=1$. Otherwise, using a sorting
algorithm and algorithm \ref{oper8}, we compute
\[
\md(u_1,v_1), \quad \lm(u_1,v_1),
\]
move to the pair of elements  $u'$, $v'$  and make the same
 calculations with the second chain (see
Section~\ref{sec:DIV} for details). If the length of words $u$ and
$v$ is bounded from above by $n$, we have to repeat these
procedures no more than $n$ times. Denote these procedures $\mathcal{GD}$
and $\mathcal{LM}$ and their complexity functions ${\rm GD}(n)$ and ${\rm LM}(n)$.
Since every step of these algorithms requires $4$ elementary
operation, we established the following fact.

\begin{prop}
\[
{\rm GD}(n) \leqslant 4n, \quad {\rm LN}(n) \leqslant 4n.
\]
\end{prop}

Let $Y \subset X$ be a proper subset of $X$. Given a geodesic word
$w$, denote by $\md_Y(w)$ the maximal left divisor of $w$ subject
to the condition that all its letters belong to $Y$. Let ${\rm GD}_Y$
denotes the time complexity of the algorithm for finding
$\md_Y(w)$.

 This is
how the algorithm works. First we compute, by means of the
algorithm for chain decomposition,  the first chain $w_1$
of $w$ and the decomposition $ w = w_1 \cdot w'$. According to
Section~\ref{sec:complexity-E'}, we need for that at most $3r$
elementary operations. Next we find the maximal $Y$-divisor $u_1$
of the abelian word $w_1$, as defined in Proposition
\ref{prop:24}, with the help of elementary operation \ref{oper8}.
Similarly, we find maximal $Y$-divisors $u_2,u_3,\dots$ of $w_2,
w_3, \dots$. The desired divisor is $u_1u_2\cdots u_k$.

If we denote, as usual, $n = l(w)$, then we come to the following
result.

\begin{prop}
\[
{\rm GD}_Y(n) \leqslant (3r+1)n.
\]
\end{prop}

\label{sec58}
This also gives us the complexity of the algorithm
$\pP_{x,Y}$ which finds the left divisor $p =
p_x(w)$ of the word $w$ maximal subject to condition that all
letters in $p$ commute with $x$.

\begin{prop}
\[
P_{x,Y}(n) \leqslant (3|Y|+1) n.
\]
\end{prop}

\subsection{The complexity of computing cyclically reduced forms}
 The algorithm $\mathcal{CR}$ takes as input a
geodesic word $w$ and outputs a cyclically reduced word $\cred(w)$
which is conjugate to $w$ and is produced as follows. Denote by
$\overleftarrow{w}$ the reverse of the word $w$ (that is, it is
written by the same letters in the reverse order). Let $w_1$ and
$\overleftarrow{w}_{1}$ be the first chains in the chain
decomposition of words $w$ and $\overleftarrow{w}$,
correspondingly. Compute
\[
\az(w_{1}) \cap \az(\overleftarrow{w}_{1}) = \{z_1,\dots,
z_{k_1}\},
\]
which is elementary operation \ref{oper10}. The total number of
elementary operations in this step of the algorithm is at most
$6r+1$. If the intersection is empty, then $\cred(w) = w$ and the
algorithm stops. Otherwise we apply $k_1$ operations of cyclic
reduction (see the fourth clause of Lemma \ref{noid}),
transforming $w$ to a geodesic word $w'$. The we repeat the same
procedure for $w'$, etc. Obviously, this process requires at most
$[n/2]$ steps. Collecting the time complexities of all stages of
reduction, we see that
\begin{prop}
\[
{\rm CR}(n) \leqslant \frac{3r}{2}n^2 +O(n).
\]
\end{prop}

\subsection{The complexity of block decomposition and finding the exhausted
form} Denote by $\mathcal{B}$ the algorithm of block decomposition which
rewrite a geodesic word $w$ as the product of blocks:
\[
w = w^{(1)}w^{(2)}\cdots w^{(k)}.\] The algorithm works as
follows.
\begin{alg} \
\begin{itemize}
\item[1\1] Find connected components $C_1,\dots,C_k$ of the graph
$\Gamma(w)$ (this is an elementary operation \ref{oper11}).
\item[2\1] Initialise $w^{(1)} = \emptyset,\dots, w^{(k)} = \emptyset$.
\item[3\1] For $i = 1,\dots, k$ do
\begin{itemize}
\item[] If $w[i] \in C_j$ then $w^{(j)} \leftarrow w^{(j)} \circ w[i]$.
\end{itemize}
\item[4\1] Return $w^{(1)},\dots, w^{(k)}$.
\end{itemize}
\end{alg}

\begin{prop}
$
{\rm B}(n) \leqslant O(r^2) + O(n).
$
\end{prop}

Now we want to estimate the complexity of the algorithm for
finding the exhausted form of a cyclically reduced word $w$. We
impose on $w$ the extra condition that $w$ is a block and also
$\az(w) = X$.

This algorithm $\mathcal{E}_Y$ was described in details in the proof of
Proposition \ref{prop:25}. The number of elementary operations it
requires does not depend on the length of $w$ and is bounded above
by a constant $|Y|$, $Y \subset X$.

\begin{prop} \qquad ${\rm E}_Y(n ) \leqslant |Y|.$
\end{prop}



\subsection{The complexity of the algorithm for the conjugacy
problem} In our usual notation, let $u$ and $v$ be two words in
the alphabet $X$. The conjugacy algorithm \cC, when applied to
words $u$ and $v$, decides, whether the elements of $G = \left<
X\right>$ represented by the words $u$ and $v$ conjugate or not.
The algorithm was described in Section~\ref{ss:scemeforconprob},
here we only discuss its complexity. We shall prove two closely
related bounds for its complexity. In one case, we estimate the
function ${\rm C}(n,r)$ which gives the number of elementary operations
\ref{oper1}--\ref{oper11} needed for deciding, in the worst case
scenario, the conjugacy of two words of length at most $n$.

Our second estimate is motivated by the following observation. The
majority of our algorithms are algorithms of dynamic
 sorting, based on comparing pairs of
letters. Therefore for the second estimate we limit the list of
elementary operations to just three:  cancellation (\ref{oper1}),
transposition of a letter and an admissible interval (\ref{oper2})
and query \ref{oper3}). All other elementary operations have
linear complexity with respect to these three basic operations.
This observation follows from 
the simple fact that, in the free
abelian group, all our algorithms work in linear time modulo these
three basic operation.
 Therefore it is natural to produce an estimate of the complexity
 of the conjugacy algorithm \cC\ in terms of the three basic
 operations.

 \subsubsection{The first estimate}
 \begin{prop} \qquad ${\rm C}(n,r) \leqslant O(r^2)\cdot O(n^2).$
\label{prop:59}
\end{prop}

\begin{proof}   The algorithm starts with
 reducing the input elements $u$ and $v$ to normal form $\n(u)$
 and $\n(v)$, which requires quadratic time in the worst case. We
 do not repeat this operation at later steps of the algorithm.
 Then we have to find block decompositions of the words $\n(u)$
 and $\n(v)$, which can be done in linear time. If the number of
 blocks is at least two, then the algorithm either gives the
 negative answer, or reduces the problem to groups of smaller
 rank. By induction, we can assume that  the normal words
 $u = \n(u)$ and $v = \n(v)$ each has only one block.
Now the algorithm finds cyclically reduced forms $\cred(u)$ and
$\cred(v)$. These operation have to be called just once, and they
can be done in quadratic time at worst. We replace $u \leftarrow
\cred(u)$ and $v \leftarrow \cred(v)$. If one of $\az(u)$ or
$\az(v)$ is a proper subset of $X$, then the algorithm either
gives negative answer, or reduces the problem to a group of
smaller rank. Therefore we assume that $\az(u) = \az(v) = X$ and
rewrite $u$ and $v$ in normal form with respect to the letter
$x_1$. This can be done in linear time modulo the timing of
elementary operations. After that the algorithm computes at most
$n$ $i$-cyclic permutations of elements in $u$, producing words
$u_1,\dots, u_k$, $k < n$, and computes exhausted forms
$e(u_1),\dots, e(u_k), e(v)$, which requires at most linear number
of operations in total, since every exhausted form can be found in
constant number of elementary operations. Finally, we have to
check $k$ equalities $e(u_1) = e(v), \dots, e(u_k) = e(v)$, which
again requires quadratic number of elementary operations. The
elements $u$ and $v$ are conjugate if and only if $e(u_i) = e(v)$
for some $i$.
\end{proof}

 \subsubsection{The second estimate} Let $\bar{\rm C}(n,r)$ be the
 number of basic operations (\ref{oper1}), (\ref{oper2}) and (\ref{oper3}) needed
 for executing the conjugacy algorithm on words of length at most
 $n$. Since elementary operations can be reduced to linearly
 bounded number of  basic operations,
 we come to the following result.
 \begin{thm} \label{thm:59} \label{complexityofconprob}
$
\bar{\rm C}(n,r) \leqslant O(r^2n^3).
$
 \end{thm}

\subsection{Complexity of the conjugation problem restricted to
cyclically reduced words} It is easy to observe that, of the
algorithms introduced in the present paper, the most complicated
from the point of view of the worst-case complexity are the
following:
\begin{itemize}
    \item the algorithm for computation of the normal form, and
    \item the algorithm for computation of the cyclically reduced
    form.
\end{itemize}
As a rule, higher level algorithms (for solving the conjugacy
problem, for instance) require very few applications of these
complicated algorithms for almost all input words. We expect, by
analogy with the work~\cite{amalgam}, that the average time
complexity of algorithms for computation of normal and cyclically
reduced forms is much lower than their worst case time complexity.

Therefore it makes sense to record the estimates for the
complexity of the algorithm for conjugacy problem with inputs
restricted to cyclically reduced words.  Analysis of the proofs of
Proposition \ref{prop:59} and Theorem
\ref{thm:59} easily leads to the following result.

\begin{thm} When restricted to cyclically reduced input words, the
algorithm \cC\  for conjugacy problem requires at most $O(rn)$
elementary operations or $O(rn^2)$ basic operations.
\end{thm}

\textbf{Acknowledgements.} The authors are grateful to Alexandre Borovik whose work on the paper far exceeded the usual editor's duties.

\section{Appendix}

Nicholas Touikan has drawn our attention to the fact that the proof
of Proposition 4.1 is incorrect. This proposition, however, is never
used in proofs of any other results of the paper and was written for
the sake of completeness of the exposition. Furthermore, Proposition
4.1 is correct (including its proof) if instead of the ShortLex
ordering on the free group $F\langle X\rangle$ (where $X$ is an
alphabet) we use yet another well ordering. Below we write out this
well ordering and the corresponding set of Knuth-Bendix rules. We
refer the reader to \cite{Epstein} for basic results on the
Knuth-Bendix procedure and use the notation introduced in the paper.

Let  $X=\{x_1, \ldots, x_n\}$. We introduce a well ordering on the
free monoid $M(X\bigcup X^{-1})$ which induces a well ordering of
the free group $F\langle X\rangle$, denote it $\prec$, as follows.
\begin{defnn}
Let $f$ and $g$ be two words from $M(X\bigcup X^{-1})$.
\[
f=u_1x_1^{\epsilon_1}\ldots u_k x_1^{\epsilon_k} u_{k+1} \quad
g=u'_1x_1^{\epsilon'_1}\ldots u'_{k'} x_1^{\epsilon'_{k'}} u'_{k'+1}
\]
where $\epsilon_i,\epsilon'_{i}\in \mathbb{Z}$, $u_j,u'_j\in M(X
\backslash \{ x_1 \} \bigcup X^{-1} \backslash \{ x_1^{-1} \})$ The
word $f$ is less than the word $g$, write $f \prec g$, if
\begin{itemize}
    \item $|f|<|g|$;
    \item $|f|=|g|$, $k<k'$;
    \item $|f|=|g|$, $k=k'$, $(\epsilon_1, \ldots, \epsilon_k)$
    precedes $(\epsilon'_1, \ldots, \epsilon'_k)$
    lexicographically;
    \item $|f|=|g|$, $k=k'$, $(\epsilon_1, \ldots, \epsilon_k)=(\epsilon'_1, \ldots,
    \epsilon'_k)$ and $u_1\prec u'_1$ or $u_1=u'_1$ and $u_2\prec u'_2$ and so on.
    Roughly saying that is $(u_1, \ldots ,u_{k+1})$ precedes $(u'_1, \ldots
    ,u'_{k+1})$ lexicographically. Note that here we assume order $\prec$ to be defined
    for $M(X\backslash \{ x_1 \} \bigcup X^{-1} \backslash \{ x_1^{-1} \})$.
\end{itemize}
\end{defnn}

One can prove the following

\begin{propn}
Relation $ \prec $ is a linear ordering of $F \langle X \rangle$.
Moreover the order $\prec$ is a well-ordering.
\end{propn}

Let $w^ \circ  $ be the $\prec$-minimal representative in the class
of all words $[w]$ equal to $w$ in $F\langle X \rangle$.

Define Knuth-Bendix rules on the set $F \langle X \rangle$ by the
following recurrent relation:
\begin{gather} \label{KBn}
\begin{split}
KB_n &= KB_{n - 1} (x_2 ,\ldots,x_n ) \cup \{  x_i^\eta
wx_1^\varepsilon \to w x_1^\varepsilon x_i^\eta ,i \ne 1; x_1 ,x_i
\notin \alpha (w);\\
&[x_i ,\alpha (wx_1 )] = 1;\varepsilon ,\eta \in \{ 1, - 1\} \}
\cup \{ x_1^\varepsilon  x_1^{ - \varepsilon } \to 1,\varepsilon \in
\{ 1, - 1\} \},
\end{split}
\end{gather}
where $w$ is written in the normal form.

\begin{propn}
The set $KB_n $ is a complete set of rules for the group $\GG_r $.
\end{propn}
\begin{proof} We use the $k$-completeness criterion from \cite{Epstein}. We
shall verify the two clauses of the criterion. We first treat the
case when two left sides of the rules overlap.

There are four cases to consider:
\begin{enumerate}
    \item \label{1} $x_{i_1 }^{\varepsilon _1 } w_1 x_{i_2 }^{\varepsilon _2 } w_2
x_{i_3 }^{\varepsilon _3 } $, where $i_1  > i_2 $, $x_{i_1 } ,x_{i_2
} \notin \alpha (w_1 )$, $[x_{i_1 } ,\alpha (w_1 x_{i_2
}^{\varepsilon _2 } )] = 1$, $i_2  > i_3 $, $x_{i_2 } ,x_{i_3 }
\notin \alpha (w_2 )$, $[x_{i_2 } ,\alpha (w_2 x_{i_3 }^{\varepsilon
_3 } )] = 1$, $\varepsilon _1 ,\varepsilon _2 ,\varepsilon _3  \in
\{ 1, - 1\} $;
    \item  \label{2} $x_{i_1 }^{\varepsilon _1 } w_1 x_{i_2
}^{\varepsilon _2 } w_2 x_{i_3 }^{\varepsilon _3 } w_3 x_{i_4
}^{\varepsilon _4 } $, where $i_1  > i_3 $, $x_{i_1 } ,x_{i_3 }
\notin \alpha (w_1 x_{i_2 }^{\varepsilon _2 } w_2 )$, $[x_{i_1 }
,\alpha (w_1 x_{i_2 }^{\varepsilon _2 } w_2 x_{i_3 }^{\varepsilon _3
} )] = 1$, $i_2  > i_4 $, $x_{i_2 } ,x_{i_4 }  \notin \alpha (w_2
x_{i_3 }^{\varepsilon _3 } w_3 )$, $[x_{i_2 } ,\alpha (w_2 x_{i_3
}^{\varepsilon _3 } w_3 x_{i_4 }^{\varepsilon _4 } )] = 1$,
$\varepsilon _1 ,\varepsilon _2 ,\varepsilon _3 ,\varepsilon _4  \in
\{ 1, - 1\} $;
    \item \label{3}
   $x_i^{ - \varepsilon } x_i^\varepsilon  wx_j^\eta $, where $j <
i$, $x_i ,x_j \notin \alpha (w_1 )$, $[x_i ,wx_j^\eta ] = 1$,
$\varepsilon ,\eta \in \{ 1, - 1\} $;
    \item  \label{4} $x_i^\varepsilon wx_j^\eta
x_j^{ - \eta } $, where $j < i$, $x_i ,x_j  \notin \alpha (w_1 )$,
$[x_i ,wx_j^\eta  ] = 1$, $\varepsilon ,\eta  \in \{ 1, - 1\} $.
\end{enumerate}

We use induction on the number of generators to verify case \ref{1}.

\[
\begin{array}{ccccccc}
                                                                                            &             \\
                                                                                            & \nearrow    \\
{x_{i_1 }^{\varepsilon _1 } w_1 x_{i_2 }^{\varepsilon _2 } w_2 x_{i_3 }^{\varepsilon _3 } } &  \\
                                                                                            & \searrow \\
                                                                                            &            \\
\end{array}
\begin{CD}
{w_1 x_{i_2 }^{\varepsilon _2 } x_{i_1 }^{\varepsilon _1 } w_2
x_{i_3 }^{\varepsilon _3 } } @>>> {w_1 x_{i_1 }^{\varepsilon _1 }
w_2 x_{i_3 }^{\varepsilon _3 } x_{i_2 }^{\varepsilon _2 } }
\\
 @. @VV{\tau}V  \\
{x_{i_1 }^{\varepsilon _1 } w_1 w_2 x_{i_3 }^{\varepsilon _3 }
x_{i_2 }^{\varepsilon _2 } } @>{\sigma}>> t
\end{CD}
\]

Transformations $\tau$ and $\sigma$ fall under the inductive
assumption, for $ x_{i_2 }  \notin \alpha (x_{i_1 }^{\varepsilon _1
} w_1 ) $, and $ x_{i_1 }^{\varepsilon _1 } w_1 $ equals $ w_1
x_{i_1 }^{\varepsilon _1 }$ in $\GG_r$.

In case (\ref{2}) write
\[
\begin{array}{*{20}c}
   {x_{i_1 }^{\varepsilon _1 } w_1 x_{i_2 }^{\varepsilon _2 } w_2 x_{i_3 }^{\varepsilon _3 } w_3 x_{i_4 }^{\varepsilon _4 } }  \\
    \downarrow   \\
   {w_1 x_{i_2 }^{\varepsilon _2 } w_2 x_{i_3 }^{\varepsilon _3 } x_{i_1 }^{\varepsilon _1 } w_3 x_{i_4 }^{\varepsilon _4 } }  \\
 \end{array}
 \begin{array}{*{20}c}
    \to   \\
   {}  \\
    \to   \\
 \end{array}
 \begin{array}{*{20}c}
   {x_{i_1 }^{\varepsilon _1 } w_1 w_2 x_{i_3 }^{\varepsilon _3 } w_3 x_{i_4 }^{\varepsilon _4 } x_{i_2 }^{\varepsilon _2 } }  \\
    \downarrow   \\
   {w_1 w_2 x_{i_3 }^{\varepsilon _3 } x_{i_1 }^{\varepsilon _1 } w_3 x_{i_4 }^{\varepsilon _4 } x_{i_2 }^{\varepsilon _2 } }  \\
 \end{array}
 \]

Consider the third case
\[\begin{array}{*{20}c}
   {x_i^{ - \varepsilon } x_i^\varepsilon  wx_j^\eta  }  \\
    \downarrow   \\
   {wx_j^\eta  }  \\
 \end{array}
 \begin{array}{*{20}c}
    \to   \\
   {}  \\
    \leftarrow   \\
 \end{array}
 \begin{array}{*{20}c}
   {x_i^{ - \varepsilon } wx_j^\eta  x_i^\varepsilon  }  \\
    \downarrow   \\
   {wx_j^\eta  x_i^{ - \varepsilon } x_i^\varepsilon  }  \\
 \end{array}
 \]

Similarly in the fourth case we get
\[
\begin{array}{*{20}c}
   {x_i^\varepsilon  wx_j^\eta  x_j^{ - \eta } }  \\
    \downarrow   \\
   {wx_i^\varepsilon  }  \\
\end{array}
\begin{array}{*{20}c}
    \to   \\
   {}  \\
    \leftarrow   \\
 \end{array}
 \begin{array}{*{20}c}
   {wx_j^\eta  x_i^\varepsilon  x_j^{ - \eta } }  \\
    \downarrow   \\
   {wx_j^\eta  x_j^{ - \eta } x_i^\varepsilon  }  \\
 \end{array}
\]

We now verify the second clause of the $k$-completeness criterion,
i.e. the case when one of the left sides of a rule is a subword of
another one. There are two cases to consider.

\begin{itemize}
    \item $x_{i_1 }^{\varepsilon _1 } w_1 x_{i_2 }^{\varepsilon _2 } w_2
x_{i_3 }^{\varepsilon _3 } $, where $i_1  > i_3 $, $x_{i_1 } ,x_{i_3
} \notin \alpha (w_1 x_{i_2 }^{\varepsilon _2 } w_2 )$, $[x_{i_1 }
,\alpha (w_1 x_{i_2 }^{\varepsilon _2 } w_2 x_{i_3 }^{\varepsilon _3
} )] = 1$, $i_2  > i_3 $, $x_{i_2 } ,x_{i_3 } \notin \alpha (w_2 )$,
$[x_{i_2 } ,\alpha (w_2 x_{i_3 }^{\varepsilon _3 } )] = 1$,
$\varepsilon _1 ,\varepsilon _2 ,\varepsilon _3  \in \{ 1, - 1\} $;

\item  $x_{i_1 }^{\varepsilon _1 } w_1 x_{i_2 }^{\varepsilon _2 }
w_2 x_{i_3 }^{\varepsilon _3 } $, where $i_1  > i_3 $, $x_{i_1 }
,x_{i_3 } \notin \alpha (w_1 x_{i_2 }^{\varepsilon _2 } w_2 )$,
$[x_{i_1 } ,\alpha (w_1 x_{i_2 }^{\varepsilon _2 } w_2 x_{i_3
}^{\varepsilon _3 } )] = 1$, $i_1  > i_2 $, $x_{i_1 } ,x_{i_2 }
\notin \alpha (w_1 )$, $[x_{i_1 } ,\alpha (w_1 x_{i_2 }^{\varepsilon
_2 } )] = 1$, $\varepsilon _1 ,\varepsilon _2 ,\varepsilon _3  \in
\{ 1, - 1\} $;
\end{itemize}

In the first case we have
\[
\begin{array}{*{20}c}
   {x_{i_1 }^{\varepsilon _1 } w_1 x_{i_2 }^{\varepsilon _2 } w_2 x_{i_3 }^{\varepsilon _3 } }  \\
    \downarrow   \\
   {x_{i_1 }^{\varepsilon _1 } w_1 w_2 x_{i_3 }^{\varepsilon _3 } x_{i_2 }^{\varepsilon _2 } }  \\
 \end{array}
 \begin{array}{*{20}c}
    \to   \\
   {}  \\
    \to   \\
 \end{array}
 \begin{array}{*{20}c}
   {w_1 x_{i_2 }^{\varepsilon _2 } w_2 x_{i_3 }^{\varepsilon _3 } x_{i_1 }^{\varepsilon _1 } }  \\
    \downarrow   \\
   {w_1 w_2 x_{i_3 }^{\varepsilon _3 } x_{i_1 }^{\varepsilon _1 } x_{i_2 }^{\varepsilon _2 } }  \\
 \end{array}
\]

Consider the second case.
\[
\begin{array}{*{20}c}
   {x_{i_1 }^{\varepsilon _1 } w_1 x_{i_2 }^{\varepsilon _2 } w_2 x_{i_3 }^{\varepsilon _3 } }  \\
    \downarrow   \\
   {w_1 x_{i_2 }^{\varepsilon _2 } w_2 x_{i_3 }^{\varepsilon _3 } x_{i_1 }^{\varepsilon _1 } }  \\
 \end{array}
 \begin{array}{*{20}c}
    \to   \\
   {}  \\
    =   \\
 \end{array}
 \begin{array}{*{20}c}
   {w_1 x_{i_2 }^{\varepsilon _2 } x_{i_1 }^{\varepsilon _1 } w_2 x_{i_3 }^{\varepsilon _3 } }  \\
    \downarrow   \\
   {w_1 x_{i_2 }^{\varepsilon _2 } w_2 x_{i_3 }^{\varepsilon _3 } x_{i_1 }^{\varepsilon _1 } }  \\
 \end{array}
\]

Therefore we checked all the conditions of the $k$-completeness
criterion. \end{proof}

Moreover the ordering $\prec$ has the following interesting
property.

\begin{propn}
For any element $w \in \GG_r $, the normal form $w^ \circ $ and the
$HNN$-normal form $n(w)$ coincide.
\end{propn}
\begin{proof} We use induction on $r$ and show that any word written in the
normal form $w^ \circ $ is in the normal form $n(w)$. Let
\[
n(w) = s_0 x_1^{\alpha _1 } s_1 x_1^{\alpha _2 } s_2 \ldots s_{k -
1} x_1^{\alpha _k } v, \quad w^ \circ = u_1 x_1^{\varepsilon _1 }
\ldots u_l x_1^{\varepsilon _l } u_{l + 1}.
\]
If we treat $n(w)$ and $w^\circ$ as elements of an $HNN$-extension
we get $k = l$, $\alpha _i = \varepsilon _i $ and $s_i $, $u_i $
represent the same coset of $A$ in $\GG_{r - 1} $.

Consider the word $n(w)$. Suppose that  $n(w)\not \simeq w^ \circ $.
Then $n(w)$ is not the $\prec$ minimal word in $[w]$. Thus, since
(\ref{KBn}) is a complete set of rules for the order $\prec$ there
exists an interval $n(w)[l,m]$ of $n(w)$ which coincides with a left
part of a rule $\rho$ from (\ref{KBn}). By the inductive assumption,
the rule $\rho$ does not lie in the set $KB_{r-1}(x_2,\dots,x_r)$,
more precisely for $s_i$ and $v$ in $n(w)$ holds:
\[
s_i^ \circ \simeq n(s_i ), v^ \circ \simeq n(v).
\]
Clearly, the rule $\rho$ can not be a cancellation. Consequently,
 $\rho$ has the form:
\begin{gather} \notag
\begin{split}
\rho=&(x_i^\eta w x_1^\varepsilon \to w x_1^\varepsilon x_i^\eta
),\hbox{ where } i \ne 1;\\
& x_1 ,x_i \notin \alpha (w);\ [x_i ,\alpha (wx_1 )] = 1;\
\varepsilon ,\eta  \in \{ 1, - 1\},
\end{split}
\end{gather}
i.e. the letter $x_i$ from some representative $s$ can be moved to
the rightmost position by the means of commutativity relations of
$G_r$. This derives a contradiction with the condition for a word to
be in the HNN-normal form.
\end{proof}

\end{document}